\documentclass[12pt]{article}
\usepackage{geometry}
\usepackage[T1]{fontenc}     
\usepackage{lmodern}         
\usepackage{amssymb}
\newcounter{bla}

\usepackage{amsmath}
\usepackage[numbers]{natbib}  
\usepackage[pdftex]{graphicx}
\usepackage{calc}
\usepackage{multicol}
\usepackage{multirow}
\usepackage{mathtools}
\usepackage{siunitx}
\usepackage{booktabs}
\usepackage{rotating}
\usepackage{wrapfig}
\usepackage[colorlinks]{hyperref}
\hypersetup{
	colorlinks,
	citecolor=magenta,
	linkcolor=cyan,
	filecolor=blue,      
	urlcolor=red,
	pdftitle={SMLMC},
	pdfpagemode=FullScreen,
}
\usepackage{etoolbox}
\usepackage{tikz}
\usepackage{url}
\usepackage{cleveref}
\usepackage[export]{adjustbox}

\newtheorem{proposition}{Proposition}[section]  
\usepackage[shortlabels]{enumitem}
\usepackage{float}
\usepackage[font=normalsize]{caption}
\usepackage{subcaption} 
\usepackage{appendix}  

\AtBeginEnvironment{thebibliography}{\sloppy}

\usepackage[dvipsnames]{xcolor}

\newcommand{\rvsg}[1]{\textcolor{black}{#1}}
\newcommand{\rvsb}[1]{\textcolor{black}{#1}}

\newcommand{\rvst}[1]{\textcolor{black}{#1}}
\newcommand{\rvstt}[1]{\textcolor{black}{#1}}
\newcommand{\rvso}[1]{\textcolor{black}{#1}}
\usepackage{pgfplots} 
\usepackage{pgfgantt}
\usepackage{pdflscape}
\pgfplotsset{compat=newest} 
\pgfplotsset{plot coordinates/math parser=false}

\definecolor{myblue}{rgb}{0.2, 0.4, 0.9}
\newcommand{\rvsn}[1]{\textcolor{black}{#1}}

\newcommand{\mrm}{\mathrm}     

\newcommand{\ex}[1]{\times 10^{#1}}

\newcommand{\vek}[1]{\mathchoice{\displaystyle\boldsymbol#1}
{\textstyle\boldsymbol#1}{\scriptstyle\boldsymbol#1}
{\scriptscriptstyle\boldsymbol#1}}

\usepackage[nopostdot,nonumberlist,acronym,toc,section]{glossaries}
\glsdisablehyper 
\newacronym{bmd}{BMD}{bone mineral density}
\newacronym{qct}{QCT}{quantitative computed tomography}
\newacronym{ct}{CT}{computed tomography}
\newacronym{hu}{HU}{Hounsfield units}
\newacronym{dxa}{DXA}{dual energy X-ray absorptiometry}
\newacronym{fea}{FEA}{finite element analysis}
\newacronym{fem}{FEM}{finite element method}
\newacronym{mc}{MC}{Monte Carlo method}
\newacronym{mpp}{MPP}{most probable point}
\newacronym{amv}{AMV}{advanced mean-value method}
\newacronym{form}{FORM}{first order reliability method}
\newacronym{sorm}{SORM}{second order reliability method}
\newacronym{rsms}{RSMs}{response surface methods}
\newacronym{pce}{PCE}{polynomial chaos expansion}
\newacronym{bc}{BC}{boundary condition}
\newacronym{pdf}{PDF}{probability density function}
\newacronym{mse}{MSE}{mean square error}
\newacronym{rmse}{RMSE}{root mean square error}
\newacronym{lhs}{LHS}{Latin hypercube sampling}
\newacronym{ism}{ISM}{importance sampling method}
\newacronym{kle}{KLE}{Kosambi-Karhunen-Lo\`{e}ve expansion}
\newacronym{spde}{SPDE}{stochastic partial differential equation}
\newacronym{pde}{PDE}{partial differential equation}
\newacronym{mlmc}{MLMC}{multilevel Monte Carlo method}
\newacronym{dof}{DOF}{degrees of freedom}
\newacronym{iso}{iso}{isotropy}
\newacronym{ortho}{ortho}{orthotropy}
\newacronym{scl}{scl}{scaling}
\newacronym{dir}{dir}{direction}
\newacronym{uq}{UQ}{Uncertanity quantification}
\newacronym{pl}{PL}{Plaston}

\newcommand{\spatialdomain}{\mathcal{G}}

\newcommand{\Euclideanspace}{\mathbb{R}^d}
\newcommand{\straintensor}{\vek{\varepsilon}(\vek{x})}
\newcommand{\dispvector}{\vek{u}}

\newcommand{\dispvectorApprox}{\vek{u}_h}
\newcommand{\stresstensor}{\vek{\sigma}(\vek{x})}
\newcommand{\volforcevector}{\vek{f}(\vek{x})}

\newcommand{\elasttensor}{\vek{C}}
\newcommand{\elastmatrix}{C}

\newcommand{\tracvector}{\vek{t}(\vek{x})}

\newcommand{\solspace}{\mathcal{U}}
\newcommand{\finitespace}{\mathcal{U}_h}

\newcommand{\elesize}{h}

\newcommand{\convorder}{\alpha}
\newcommand{\convbeta}{\beta}
\newcommand{\convgamma}{\gamma}
\newcommand{\convorderv}{\alpha}
\newcommand{\convbetav}{\beta}

\newcommand{\calpha}{c_{\alpha}}
\newcommand{\cbeta}{c_{\beta}}
\newcommand{\cgamma}{c_{\gamma}}
\newcommand{\ccost}{c}
\newcommand{\samplespace}{\varOmega}
\newcommand{\event}{\omega}

\newcommand{\probability}{\mathbb{P}}

\newcommand{\erf}{\text{erf}}

\newcommand{\samplesize}{N}

\newcommand{\corrlength}{l_c}

\newcommand{\coeffdispC}{\delta_{\elastmatrix}}
\newcommand{\coeffdispT}{\delta_{{T}}}

\newcommand{\CRF}{C(\vek{x},\omega)}

\newcommand{\stogermnoX}{{T}(\vek{x}, \event)}

\newcommand{\uppermean}{{U}}
\newcommand{\upperT}{{V}(\vek{x},\omega)}

\newcommand{\identitymatrix}{{I}}

\newcommand{\normconst}{N_c}
\newcommand{\gammaCDF}{F_{\Gamma_\alpha}}

\newcommand{\meshlevel}{l}

\newcommand{\numofmesh}{{L}}

\newcommand{\optSam}{N_l}

\newcommand{\costlevel}{\mathcal{C}_\meshlevel}
\newcommand{\costfunction}{f(\N_l)}
\newcommand{\Lagrange}{\tau}

\newcommand{\GammaRF}{\Gamma_{\alpha_j}(\vek{x},\event)}

\newcommand{\UhlNO}{u_{\h_\meshlevel, 50}^{(t)}}

\newcommand{\uu}{u} 
\newcommand{\h}{h}
\newcommand{\N}{N}
\newcommand{\uh}{\uu_\h}

\newcommand{\uhNt}{\uu_{\h,\N}^{(t)}}
\newcommand{\U}{\mathcal{U}}
\newcommand{\mumu}{\mu} 
\newcommand{\muMC}{\widehat{\mumu}^{\mathrm{MC}}} 
\newcommand{\var}{\mathbb{V}\text{ar}}

\newcommand{\htwoMC}{\widehat{\mathrm{h}}_2^{\mathrm{MC}}}
\newcommand{\htwoML}{\widehat{\mathrm{h}}_2^{\mathrm{ML}}}
\newcommand{\muML}{\widehat{\mumu}^{\mathrm{ML}}} 
\newcommand{\Vtwo}{{\mathbb{V}}_2}
\newcommand{\VtwoMC}{\widehat{\mathbb{V}}_2^{\mathrm{MC}}}
\newcommand{\Nl}{\N_\meshlevel}
\newcommand{\LL}{L} 
\newcommand{\uhO}{\uu_{\h_0}}
\newcommand{\uhl}{\uu_{\h_\meshlevel}}
\newcommand{\uhll}{\uu_{\h_{\meshlevel-1}}}
\newcommand{\uhL}{\uu_{\h_\LL}}
\newcommand{\uhONO}{\uu_{\h_0,\N_0}}
\newcommand{\uhONOt}{\uu_{\h_0, N_0}^{(t)}}
\newcommand{\uhlNl}{\uu_{\h_\meshlevel,\Nl}}

\newcommand{\uhllNl}{\uu_{\h_{\meshlevel-1},\Nl}}
\newcommand{\uhLNL}{\uu_{\h_\LL,\left\lbrace \Nl \right\rbrace}}
\newcommand{\uhLNLt}{\uu_{\h_\LL,\left\lbrace \Nl \right\rbrace}^{(t)}}

\newcommand{\YhlNl}{Y_{\meshlevel}}
\newcommand{\Yl}{Y_{\meshlevel}}

\newcommand{\Vltwo}{\mathbb{V}_{\meshlevel,2}}
\newcommand{\VltwoMC}{\widehat{\mathbb{V}}_{\meshlevel,2}^{\mathrm{MC}}}
\newcommand{\Zl}{{Z}_{\meshlevel}}

\newcommand{\totaldisp}{{u}_{h}^{(t)}}
\newcommand{\kurt}{\mathrm{Kurt}}
\newcommand{\qoi}{X}

\newcommand{\hrMC}{\widehat{\mathrm{h}}_\p^{\mathrm{MC}}}
\newcommand{\p}{p}

\newcommand{\s}{s}

\newcommand{\ee}{e}

\newcommand{\epsdot}{\mathring{\varepsilon}}
\newcommand{\e}{\mathring{\ee}}
\newcommand{\eMC}{\widehat{\e}}

\newcommand{\epssq}{\mathring{\epsilon}^2}
\newcommand{\epssqtwo}{\epssq/2}

\newcommand{\constA}{a}
\newcommand{\constB}{b}
\newcommand{\Norm}{\lambda}
\newcommand{\Nm}{\Norm_m}
\newcommand{\Nmhat}{\widehat{\Norm}_m}
\newcommand{\Nmhatinf}{\widehat{\Norm}_m}
\newcommand{\Nv}{\Norm_v}
\newcommand{\Nvhat}{\widehat{\Norm}_v}
\newcommand{\Nvhatinf}{\widehat{\Norm}_v}
\newcommand{\SymMatSpace}{\mrm{Sym}}
\newcommand{\PosSymMatSpace}{\SymMatSpace^+}

\newcommand{\PosSymSecRankSpaceElast}{\PosSymMatSpace(n)}
\newcommand{\SymSecRankSpace}{\SymMatSpace(d)}

\newcommand{\KLEterms}{M}

\newcommand{\feq}[1]{Eq.~(\ref{#1})} 
\newcommand{\feqs}[2]{Eqs.~(\ref{#1}) and (\ref{#2})} 
\newcommand{\feeqs}[2]{Eqs.~(\ref{#1})-(\ref{#2})} 
\newcommand{\fsec}[1]{Section~\ref{#1}}
\newcommand{\fsecs}[2]{Sections~(\ref{#1}) and (\ref{#2})} 

\newcommand{\ffig}[1]{Fig.~\ref{#1}} 
\newcommand{\ffigs}[2]{Figs.~\ref{#1} and \ref{#2}} 
\newcommand{\ftbl}[1]{Table~\ref{#1}}
\newcommand{\BRrefeq}[1]{Eq.~(\ref{#1})}
\newcommand{\prop}[1]{Proposition.~\ref{#1}} 
\newcommand{\props}[2]{Propositions.~\ref{#1} and \ref{#2}} 
\newcommand{\app}[1]{Appendix~\ref{#1}} 

\title{\rvsn{Scale-invariant Monte Carlo and multilevel Monte Carlo estimation of mean and variance: An application to simulation of linear elastic bone tissue}}
\author{Sharana Kumar Shivanand\textsuperscript{a}\thanks{Corresponding author. {E-mail:} \texttt{sharana.shivanand@kit.edu}},  Bojana Rosi\'{c}\textsuperscript{b} \\
	\small \textsuperscript{a}Institute of Scientific Computing, Technische Universit\"{a}t Braunschweig, Germany\\
	\small \textsuperscript{b}Applied Mechanics and Data Analysis, University of Twente, Netherlands}
\date{\today}

\newenvironment{keyword}{%
	\noindent \emph{Keywords:}}{\newline}

\begin{document}
\sloppy
\maketitle
    \begin{abstract}
       
      \rvsn{We propose novel scale-invariant error estimators for the Monte Carlo and multilevel Monte Carlo estimation of mean and variance. For any linear transformation of the distribution of the quantity of interest, the computation cost across \rvso{fidelity} levels is optimized using a normalized error estimate, which is not only fully dimensionless but also remains robust to variation in characteristics of the distribution.  
      We demonstrate the effectiveness of the algorithms through application to a mechanical simulation of linear elastic bone tissue, where material uncertainty incorporating both heterogeneity and random anisotropy is considered in the constitutive law.}
       

    \end{abstract}
    


%
%


\section{Introduction}\label{Introduction}

	\rvsn{
	The estimation of statistical moments for characterizing a probabilistic quantity of interest (QoI), such as the solution to a stochastic partial differential equation (SPDE) or an ordinary differential equation (SODE), is of fundamental importance in the field of uncertainty quantification (UQ) \cite{xiu2010numerical, lord_powell_shardlow_2014}. Specifically, let $\qoi$ be a random variable defined on a probability space $(\varOmega, \mathfrak{F}, \mathbb{P})$, where $\varOmega$ denotes the sample space, $\mathfrak{F}$ is a $\sigma$-algebra of measurable events, and $\mathbb{P}$ is the associated probability measure. The $\p$-th central moment of $\qoi$, for $\p \in \mathbb{N}$, is defined as
	\begin{equation}\label{eq:momentp}
		\mu_{\p}(\qoi) = \mathbb{E}\left[(\qoi - \mathbb{E}(\qoi))^{\p}\right],
	\end{equation}
	where, $\mathbb{E}(\qoi):={\mu}(\qoi)$ is the mean or the first raw moment  (not to be confused with the first central moment $\mu_1$, which equals to zero when $\p=1$). The second central moment, when $\p=2$, $\mu_2(\qoi) = \var(\qoi):=\mathbb{E}\left((\qoi-\mathbb{E}(\qoi))^2\right) $, denotes  the variance. 
	In contrast, third and fourth central moments are typically expressed in standardized form: the skewness and kurtosis, which are defined as $\alpha_3 = \mu_3 /{\sqrt{\var^3}}$ and $\alpha_4 = \mu_4 / \var^2$, respectively.
	While skewness and kurtosis are dimensionless and invariant under linear transformation of $X$, the mean and variance are inherently scale-dependent. Consequently, their estimation accuracy, typically evaluated using absolute error metrics such as the mean squared error (MSE) or root mean squared error (RMSE), also depends on the scale of the QoI. This scale dependency presents an interpretability challenge, especially in practical scenarios where comparisons across different estimators or between the same estimator applied to differently scaled QoI are required. Addressing this, normalized error estimates play a crucial role in statistics by providing a standardized measure of the accuracy of statistical estimates.}
	
	\rvsn{Normalization in statistics is a broad topic with multiple interpretations \cite{kenney1954mathematics}. One common usage refers to the standardization of observational data (or the QoI), often represented by the z-score or standard score, defined as $z = {\qoi - \mu}/{\sqrt{\var}}$. Another interpretation involves the normalization of statistical moments of the form, $\alpha_p = \mathbb{E}(z^{\p}) = \mu_p/{\sqrt{\var^p}}$, 	
	whereby the first and second standardized moments, $\alpha_1$ and $\alpha_2$--- representing the mean and variance of standardized variable $z$---are fixed at 0 and 1, respectively. 
	This implies that these moments are invariant across different distributions, offering no distinctive information. In contrast, higher-order standardized moments, such as skewness $\alpha_3$ and kurtosis $\alpha_4$, retain distribution-specific characteristics and are thus useful for distinguishing among probability distributions.
	Therefore, in this work, the primary focus is on the estimation of the mean $\mu$ and variance $\mu_2$, where normalization is applied not to the QoI or the moments, but rather to their scale-dependent absolute error estimates.} 
	In this study, we focus on random sampling-based statistical estimators, such as
	\rvst{the Monte Carlo (MC) method}; traditionally regarded as the gold standard for solving stochastic problems, due to its simplicity and resilience against the curse of dimensionality \cite{metropolis_monte_1949, Hammersley1965, fishman_monte_1996,graham_stochastic_2013}. However, despite its appeal, the practical implementation of MC often faces challenges due to slow convergence, thereby demanding substantial computational effort. To overcome these limitations, recent research has increasingly turned to variance reduction techniques, notably \rvst{the multilevel Monte Carlo (MLMC) method}. 
	A key objective of the MLMC method is to spread out the sampling strategy across a hierarchy of different fidelities (or levels) such that the number of stochastic samples drastically decays with the increment in fidelity of the model. Under the right conditions, this results in an overall reduction of the computational cost as compared to the MC approach.

	To the best of the authors’ knowledge, \rvst{MLMC} was first introduced by \cite{Heinrich2001} in the context of estimating multi-dimensional parameter-dependent integrals and is further employed by \cite{giles2008} for solving It\^{o}'s stochastic ordinary differential equations in the computational finance application. Following this, \rvst{MLMC} is extended to solve the linear elliptic PDEs describing the subsurface flow with inhomogeneous stochastic parameters \cite{Cliffe2011,Barth:2011}. \rvst{A further analysis of modelling of rough random field coefficients of elliptic PDEs on MLMC convergence} is studied in \cite{charrier_finite_2013, teckentrup_further_2013}. However, previously mentioned works focussed only on the approximation of the QoI's sample mean; therefore, lacked full characterization of the probabilistic solution. Consequently, \cite{mishra_multi-level_2012} studies the unbiased MLMC sample variance estimator, which is further analysed in \cite{bierig_convergence_2015} on a class of elliptic random obstacle problems---{nonetheless, the corresponding error estimates are defined only on worst case bounds}. To address this limitation, an alternative MLMC variance estimator based on h-statistics \cite{dwyer1937,Colin_Rose:2002} is introduced in \cite{krumscheid_quantifying_2020}, {which emphasizes on unbiased construction of the MSE in closed form, an approach also adopted in this study}. \rvst{More recently, \cite{menhorn_multilevel_2024} developed MLMC estimators for the biased standard deviation and its linear combination with the unbiased mean, which are important for optimization under uncertainty (OUU) workflows.} Though the estimation of other higher-order moments such as skewness and kurtosis or covariance structures of the QoI via MLMC is out of the scope of this paper, more literature on this can be found in \cite{bierig_estimation_2016,krumscheid_quantifying_2020, shivanand_covariance_2025}.

	In the context of assessing the convergence of the MLMC algorithm in a standardized manner, the authors in \cite{bierig_convergence_2015} theoretically prove, with specific assumptions on the deterministic solution, two error bounds for the multilevel sample mean and variance that may be directly compared. But no practical normalized error estimates are defined for interpretation of the complexities between the moments. On the other hand, relative error estimates of the MC and MLMC methods, where the total MSE of a given moment is normalized by the square of its own statistic, are considered in \cite{krumscheid_quantifying_2020}. \rvsn{However, such relative errors are not fully scale-invariant. 
	For example, in the MC estimation of the mean, any additive linear transformation of the QoI changes the ratio of absolute MSE to the squared mean in a dimension-dependent manner. Interestingly, in the MC variance estimation, one does obtain a fully dimensionless relative error under linear scaling and translation, when normalized by the square of the variance. However, the resulting error estimate becomes highly sensitive to the tail behaviour (kurtosis) of the QoI distribution.}
	To tackle this issue, in this article, we propose novel normalized mean square error (NMSE) estimates for both the MC mean and variance estimations, based on which the MLMC counterparts are derived.
	The newly introduced relative errors are statistically defined {using h-statistics}, with chosen normalizing factors that are finite and unbiased. They ensure that the total MSEs of the MC and the MLMC algorithms are fully scale-invariant under any linear transformation (scaling and addition) of QoI, and remain robust to variations in distributional characteristics. 
	Therefore, the proposed NMSEs enable easier interpretation of statistical accuracy and efficiency between the MC and MLMC algorithms for the estimation of both mean and variance, and across different scales.
    
	
    The other objective of this paper is to investigate the applicability of the scale-invariant MC and MLMC method to linear elliptic problems described by stochastic \rvso{material} parameters representing both heterogeneity and uncertain symmetries. As an example, we consider the linear-elastic material model of the human femoral bone, whose constitutive law is assumed to be uncertain. In particular, bone tissue is not only a highly heterogeneous material but also anisotropic, where the material symmetry is described as uncertain due to the lack of conclusive verification of the class of elastic symmetry it belongs to \cite{Sansalone:2016,KELLER19941159}. Hence, the entire elasticity tensor is constructed as random with a predefined elastic symmetry in the mean (orthotropy \cite{Geraldes:2014} in this case) and triclinic symmetries in each stochastic realization. Consequently, the positive-definite random elasticity matrices are modelled as matrix-valued random fields, as proposed in \cite{soize2006non,soize2008}. In this work, we restrict ourselves to material uncertainties only, assuming that the remaining model parametrization is deterministic. However, the numerical approach is general enough to be employed for other types of uncertainties as well. 
    
    The paper is organised as follows: In \fsec{UncertaintyQuantification}, we describe the problem, whereas the theoretical procedures of scale-invariant MC mean and variance estimates are elaborated in \fsec{Sec:MC}. Following this, in \fsec{Sec:MLMC}, the normalized versions of the MLMC mean and variance estimators are detailed. A deterministic and stochastic setting of the linear elastic material model, along with the stochastic material modelling, is given in \fsec{Sec:application}. In \fsec{NumericalResults}, we visualize the considered stochastic material model when implemented on a two-dimensional proximal femur, and detail the accuracy and efficiency of the normalized MC and MLMC. Finally, the conclusions are drawn in \fsec{Conclusion}.


\section{Problem description}\label{UncertaintyQuantification}

Let us consider a physical system occupying the spatial domain $\spatialdomain\subset\Euclideanspace$ in a $d$-dimensional Euclidean space, modelled by an abstract equilibrium equation:
\begin{equation}\label{abs_main_eq}
	\mathcal{A}(q(x),\uu(x))=f(x).
\end{equation}
Here, $\uu(x) \in \U$ describes the state of the system at a spatial point $x\in\spatialdomain$ lying in a Hilbert space $\U$ (for the sake of simplicity), 
$\mathcal{A}$ is a (possibly non-linear) operator modelling the physics of the system, and $f\in \U^*$, \rvst{the dual space of $\U$}, is some external influence (action/excitation/loading). Furthermore, we assume that the model depends on the parameter set $q \in \mathcal{Q}$ and that it is accompanied by appropriate boundary and/or initial conditions.
Note that, for brevity, the model in the previous equation describes the physical system only in a spatial domain $x$, whereas one may also generalise this to time-dependent processes.

The uncertainty in the previous equation may arise due to the randomness in external influence $f$,  initial or boundary conditions, geometry $\spatialdomain$, as well as the coefficient of operator $\mathcal{A},$ i.e., parameter $q$. Although the theory presented further does not depend on this choice and is general enough to cover all of the mentioned (single or combination of) cases, this article focuses on incorporating stochasticity only in coefficient $q$. In the theory of continuum solid mechanics, the parameter $q$ represents one of the very well-known physical phenomena, such as elasticity \cite{ostoja-starzewski_microstructural_2007}, which is detailed as an example in \fsec{Sec:application}. Here, we assume that $q$  
is modelled as a
random field $	q(x,\omega)$ with finite second-order moments on a probability space $(\varOmega, \mathfrak{F}, \mathbb{P})$.
Following this, \BRrefeq{abs_main_eq}
rewrites to a stochastic form:
\begin{equation}\label{basic_eq}
	\mathcal{A}(q(x,\omega),\uu(x,\omega))=f(x),
\end{equation}
which further is to be solved for $\uu(x,\omega) \in \U\otimes L_2(\varOmega, \mathfrak{F}, \mathbb{P})$. 

Often, when the response of a model is uncertain, one is interested in attaining its relevant information via its corresponding statistics. 
Analytically, the statistical moments of the solution $\uu \equiv\uu(x,\omega)$ can be represented, by rewriting \feq{abs_main_eq}, as
$\p$\textendash{th} central moment:
\begin{equation}\label{math_Exp}
	\mu_{\p}(\uu) = \mathbb{E}\left((\uu-\mathbb{E}(\uu))^{\p}\right) = \int_\samplespace (\uu-\mathbb{E}(\uu))^{\p} \probability(\text{d}\event).
\end{equation}
The objective of this study is to determine the mean and variance only after an appropriate deterministic discretization of the problem in \BRrefeq{basic_eq} is presented.

Due to spatial and stochastic dependence, the solution $\uu \equiv\uu(x,\omega)$ in \BRrefeq{basic_eq} is first discretized in a spatial domain, i.e.,~we search for the solution in a finite subspace $\U_\h\subset \U$.  After rewriting the problem in \BRrefeq{basic_eq} in a variational  form, the spatial discretization $\uh(x, \omega)\in \U_\h$---$\h$ being the discretization parameter---can take the finite element form  on a sufficiently fine spatial \rvsb{mesh} $\mathcal{T}_h$ \footnote{Other types of discretization techniques can be considered as well} \cite{zienkiewicz_finite_2013}. The expectation functional in \BRrefeq{math_Exp}  then rewrites to:
\begin{equation}\label{math_Exp_h}
	\mu_{\p}(\uh) = \mathbb{E}((\uh-\mathbb{E}(\uh))^{\p}) = \int_\samplespace (\uh-\mathbb{E}(\uh))^{\p} \probability(\text{d}\event),
\end{equation}
in which the exact solution $u$  is substituted by a semi-discretized solution $\uh \equiv \uh(x, \omega)$. 


\section{Monte Carlo method}\label{Sec:MC}

Due to the complexity of integrating 
\feq{math_Exp_h} analytically, in this paper, the focus is on the traditional sampling-based Monte Carlo (MC) method \cite{metropolis_monte_1949, Hammersley1965, fishman_monte_1996,graham_stochastic_2013}. 

\subsection{Monte Carlo estimation of mean}\label{MC}

For an unbiased estimate of the statistic $g(\omega)$, one may take a symmetric function
\begin{equation}\label{symmetric}
	\muMC(g)=\frac{1}{\N}(g(\omega_1)+g(\omega_2)+...+g(\omega_N)),
\end{equation}
meaning that the estimate does not depend on the order in which observations were taken. The existence and uniqueness of such a choice are given in \cite{halmos1946}. Under the assumption that each sample $\uh(x, \event_i)$ comes from the identical distribution as $\uh(x, \event)$ and by use of \BRrefeq{symmetric}, one may reformulate
the \acrshort{mc} estimate of the mean $\mu(\uh)$ as 
\begin{equation}
	\mathbb{E}(\uh) \approx	\muMC(\uh)  = \frac{1}{\N}\sum_{i=1}^{\N}\uh(x, {\event}_i),
\end{equation}
where $\samplesize > 1$ is the sample size of the random field $\uh(x, {\event})$ at a spatial location $x$.
Following this, the approximation error of the MC-based mean estimate $\muMC \equiv \muMC(\uh)$ compared to the exact mean reads
\begin{equation}
	{\ee}(\muMC) = \muMC -\mathbb{E}(\uu),
\end{equation}
which further can be rewritten as
\begin{equation*}\label{error}
	\ee(\muMC) = \muMC -\mathbb{E}(\uh) +\mathbb{E}(\uh) -\mathbb{E}(\uu).
\end{equation*}
Thus, the \acrfull{mse} reads
\begin{equation*}
	\begin{aligned}
		\textrm{MSE}(\muMC) :=\mathbb{E}({\ee}^2) 
		&= \mathbb{E}((\muMC -\mathbb{E}(\uh) +\mathbb{E}(\uh) -\mathbb{E}(\uu))^2) \\
		&= \mathbb{E}((\muMC -\mathbb{E}(\uh))^2)+(\mathbb{E}(\uh) -\mathbb{E}(\uu))^2+ \\
		& \qquad 2(\mathbb{E}(\uh) -\mathbb{E}(\uu))(\mathbb{E}(\muMC) -\mathbb{E}(\uh)),
	\end{aligned} 
\end{equation*}
in which $\muMC$ is assumed to be an unbiased estimator; meaning that
$ \mathbb{E}(\muMC) = \mathbb{E}(\uh) $, 
and also,
$ \mathbb{E}((\muMC -\mathbb{E}(\uh))^2) = \var(\muMC) = {\var(\uh)}/{\N} $, via the central limit theorem \cite{fischer2011history}. Therefore, the previous equation reduces to:
\begin{equation}\label{Eq:MCerror1}
	\textrm{MSE}(\muMC) = \frac{\var(\uh)}{\N}+
	(\mathbb{E}(\uh) -\mathbb{E}(\uu))^2. 
\end{equation}
The first term in the previous equation is the sampling error, which varies inversely to the sample size $\samplesize$. On the other hand, the second term describes square of the spatial discretization error, whose change is proportional to the element size $\elesize$.

\rvstt{\textbf{Scale-invariant error estimator for {MC} mean estimation}:} The MSE in \feq{Eq:MCerror1} is an absolute error estimate and, hence, scale-dependent. For example,
if each observation $\uh(x, \omega_i)$ is
linearly transformed in the form 
$	\constA\,\uh(x, {\omega}_i)+\constB, $
in which $\ \constA,\constB\in\mathbb{R}$ are constants, 
then the Monte Carlo error in \BRrefeq{Eq:MCerror1} will be affected in a square proportional manner---as 
\begin{equation}\label{Eq:scalesampleM}
	\var(\constA\,\uh+\constB)=\constA^2\,\var(\uh)
\end{equation}
and
\begin{equation}\label{Eq:scalediscM}
	(\mathbb{E}(\constA\,\uh+\constB) -\mathbb{E}(\constA\,\uu+\constB))^2 = \constA^2(\mathbb{E}(\uh) -\mathbb{E}(\uu))^2.
\end{equation}

\rvsn{This leads to an interpretability issue, which is not suitable in practical applications---for example, when transforming the units of temperature from Kelvin to Fahrenheit or displacement field from millimetre to metre. 
Moreover, an additional complication arises when comparing the convergence behaviour of different MC moments under a fixed MSE. Specifically, under a linear transformation of the form defined in the previous equations, the MSE of the MC mean estimator scales with $\constA^2$, while the MSE of the MC variance estimator scales with $\constA^4$---more information on this to follow in \fsec{MCvar}.
To address this limitation, it is necessary to define a scale-invariant version of MSE, defined in \feq{Eq:MCerror1}, by incorporating a statistical normalization factor $\Nm$, which must satisfy the following properties:}
\rvsn{\begin{enumerate}[(a)]
		\item it must be greater than zero and finite i.e., $0 < \Nm < \infty$;
		\item it should satisfy the scaling condition $\Nm(\constA\,\uh+\constB)=\constA^2\,\Nm$, {thereby achieving a complete dimensionless MSE};
		\item \rvsn{it must be chosen such that the sampling error of the scale-invariant error estimator remains robust and fully invariant with respect to variations in the properties of the distribution of solution $\uh$.}
\end{enumerate}}

\rvsn{\textbf{$\mu^2$ as a normalizer}:
It is important to highlight that the common practice of using the squared mean value {$\mu^2 = \mathbb{E}(\uh)^2$ as the standardizing quantity $\Nm$} 
satisfies the second criterion listed above—namely, achieving a dimensionless MSE—only in the case of multiplicate scale-change of the form $\mu(a\uh)^2 = a^2 \mu(\uh)^2$. However, this choice fails under additive transformations, as $\mu(\uh + b)^2 = (\mu(\uh) + b)^2$; thereby, violating the required invariance. Consequently, $\mu^2$ cannot be considered a suitable normalizing factor for constructing a truly scale-invariant error metric.} 



\rvsn{We propose using the variance of the solution, $\Nm = \var(\uh)$, as a normalizer. 
	Provided that $0 < \var(\uh) < \infty$, the factor $\Nm$ satisfies the second condition outlined above (see \feq{Eq:scalesampleM}).}
With this, the new standardized error estimate is defined as
\begin{equation}\label{new_error}
	\e(\muMC) := \frac{\ee(\muMC) }{\sqrt{\var(\uh)}}.
\end{equation}
Finally, by normalizing \BRrefeq{Eq:MCerror1}, one obtains the squared error estimate 
\begin{equation}\label{Eq:MCerror3}
	\textrm{NMSE}(\muMC) :=\mathbb{E}({\e}^2) = \frac{1}{N}+  \frac{(\mathbb{E}(\uh)-\mathbb{E}(\uu))^2}{\var(\uh)},
\end{equation}
in which the first term is the scale-invariant sampling accuracy, and the second term represents the new normalized squared discretization error. 
\rvsn{In the previous equation, the resulting NMSE remains fully dimensionless under both multiplicative and additive scale transformations. Furthermore, the normalized sampling error is only proportional to $1/\samplesize$---making it invariant to the variations in variance of solution $\uh$. In other words, for a given total NMSE and a fixed spatial resolution $\elesize$, the computational cost of the MC mean estimator $\muMC$ (see \prop{Prop:MCmean}) is solely dependent on the number of MC samples $\samplesize$, thereby fulfilling the third listed requirement.}

The computational complexity of the scale-invariant MC mean estimate is similar to the conventional procedure, as shown in \cite{Cliffe2011}, but it must be redefined with respect to the total NMSE, which is given below.
\begin{proposition}\label{Prop:MCmean}
	Let us consider $\calpha, \cgamma, \alpha$ and $\convgamma$ as positive constants, and then one may define the error bounds as follows, where
	\begin{equation*}\label{eq:errormeanmc}
		\begin{split}
			(i)& \ \text{the deterministic error decays as}\ \frac{|\mathbb{E}(\uh)-\mathbb{E}(\uu)| }{\sqrt{\Nm}} \le \calpha \h^{\alpha},\\
			(ii)& \ \text{the computational cost to determine a single realization of $\uh(x, \event)$ is given by} \\
			& \mathcal{C}(\uh)  \le \cgamma \h^{-\convgamma}.
		\end{split}
	\end{equation*}
	Then, for any $0<\mathring{\epsilon}<e^{-1}$, the Monte Carlo (MC) mean estimator $\muMC$ with $\samplesize = \mathcal{O}(\mathring{\epsilon}^{-2})$ and $\elesize = \mathcal{O}(\mathring{\epsilon}^{1/\alpha})$ satisfies the normalized mean square error $\mathrm{NMSE}(\muMC) < \epssq$. Therefore, the corresponding computational cost of MC mean estimation is
	\begin{equation*}
		\mathcal{C}(\muMC) \le  \ccost \, \mathring{\epsilon}^{-2-\convgamma/\alpha},
	\end{equation*}
	where $\ccost>0$ is a positive constant. 
\end{proposition}

\rvsn{{\textbf{Unbiased estimation of $\var$}}: In general, the true variance $\var(\uh)$ is unknown, and must be estimated. For this purpose, one may use the estimator $\widehat{\mathbb{V}}\text{ar}^{\textrm{MC}}(\uh)$, which is determined by the \acrshort{mc} procedure with $N$ random draws. However, 
	in such a case, the use of a symmetric function as in \BRrefeq{symmetric} does not lead to an unbiased estimator. To address this limitation, we make use of h-statistics for estimating the central moments $\mu_p$, which are not only unbiased but also symmetric, and possess minimal variance \cite{dwyer1937, Colin_Rose:2002}. Further details on h-statistics for univariate central moments, particularly the second-order moment (variance), are provided in \rvst{\app{appendix:hstat}}. Based on this, we rewrite \feq{new_error} by estimating $\Nm = \var(\uh)$ using the second h-statistic, denoted as $\Nmhat = \htwoMC$ (see \feq{Eq:h-stat}), resulting:
	\begin{equation}\label{new_errorMEAN}
		{\eMC(\muMC)}= \frac{\ee(\muMC)}{\sqrt{\htwoMC}}.
	\end{equation}
	Finally, one obtains a new formulation of the estimated scale-invariant MSE (from \feq{Eq:MCerror3}) in the form:
	\begin{equation}\label{Eq:MCmeanNMSE}
		\textrm{NMSE}(\muMC) = \frac{1}{N}+  \frac{(\mathbb{E}(\uh)-\mathbb{E}(\uu))^2}{\htwoMC}.
	\end{equation} }

\subsection{Monte Carlo estimation of variance}\label{MCvar} 

The previously derived estimate of variance by h-statistics in \feq{Eq:h-stat} is characterized by an approximation error given as
\begin{equation*}\label{Eq:MCvarerr1}
	{\ee}(\htwoMC)=\htwoMC-\var(\uh)+\var(\uh)-\var(\uu),
\end{equation*}
in which $\var(u)$ denotes the exact variance of the solution $u$.  
Following this, the MSE of the variance estimator reads:
\begin{equation*}\label{Eq:MCerror2}
	\begin{aligned}
		\mathrm{MSE}(\htwoMC) := \mathbb{E}({\ee}^2) 
		&= \mathbb{E}((\htwoMC-\var(\uh)+\var(\uh)-\var(\uu))^2) \\
		&= \mathbb{E}((\htwoMC-\var(\uh))^2) + (\var(\uh)-\var(\uu))^2 + \\
		& \qquad 2 (\var(\uh)-\var(\uu)) (\mathbb{E}(\htwoMC)-\var(\uh)).
	\end{aligned}
\end{equation*}
Having that $\htwoMC$ is an unbiased estimator, meaning that $\mathbb{E}(\htwoMC) = \var(\uh),$
one may further rewrite the previous equation to 
\begin{equation}\label{Eq:MCerrorvar}
	\mathrm{MSE}(\htwoMC) =\var(\htwoMC) +  ({\var}(\uh)-{\var}(\uu))^2.
\end{equation}
In the previous equation, the variance of the second-order statistic $\htwoMC$, \rvst{which also represents the normalizing constant $\Nmhat$}, is derived by \cite{E_cho_2009} and reads
\begin{equation}\label{Eq:varvar}
	\var( \htwoMC) := \frac{1}{N}\left( \mu_4(\uh) - \frac{\mu_2(\uh)^2(N-3)}{N-1}\right), 
\end{equation}
where $\mu_2(\uh)$ and $\mu_4(\uh)$ represent the second and fourth central population moments, respectively. 
For brevity, one may further rewrite the previous equation as
\begin{equation}\label{Eq:V2N}
	\var( \htwoMC) := \frac{\Vtwo(\uh)}{N},
\end{equation}
where
\begin{equation}\label{Eq:V2}
	\Vtwo(\uh) := \mu_4(\uh) - \frac{\mu_2(\uh)^2(N-3)}{N-1}.
\end{equation}
Consequently, by substituting \feq{Eq:V2N} in \feq{Eq:MCerrorvar}, one obtains
\begin{equation}\label{Eq:MCerrorvar1}
	\textrm{MSE}(\htwoMC) = \frac{\Vtwo(\uh)}{N} +  ({\var}(\uh)-{\var}(\uu))^2.
\end{equation}
Here, the first term defines the statistical error, which is inversely proportional to the number of samples $\samplesize$, and the second part represents the square of the discretization error, directly proportional to parameter $\elesize$. 

\rvstt{\textbf{Scale-invariant error estimator for {MC} variance estimation}:} Similar to MSE formulation for the MC mean estimator in \feq{Eq:MCerror1}, the error estimate for the MC variance estimator, as presented in \feq{Eq:MCerrorvar1}, is inherently scale-dependent. That is, 
under a linear transformation of the form $\constA \uh + \constB$, as discussed in \fsec{MC}, the corresponding MSE exhibits a bi-quadratic scaling behaviour. This is due to the characteristics of the central moments involved in \feq{Eq:V2}:
\begin{equation}
		\mu_4(\constA\,\uh+\constB) = \constA^4\,\mu_4(\uh) 
\end{equation}
and
\begin{equation}\label{Eq:mu2square}
	{\mu_2}^2(\constA\,\uh+\constB) = \constA^4\,{\mu_2}^2(\uh).
\end{equation}
Therefore, it is clear that
\begin{equation}\label{Eq:V2scale}
	\Vtwo(\constA\,\uh+\constB)=\constA^4\,\Vtwo(\uh)
\end{equation}
and also that
\begin{equation}
	(\var(\constA\,\uh+\constB) -\var(\constA\,\uu+\constB))^2 = \constA^4(\var(\uh) -\var(\uu))^2.
\end{equation}
Therefore, similar to the standardizing entity $\Nm$ in \fsec{MC}, the normalization entity for the MC variance estimation, denoted by $\Nv$, must satisfy the following criteria:
\rvsn{\begin{enumerate}[(a)]
		\item it must be finite and greater than zero, i.e., $0 < \Nv < \infty$;
		\item it should satisfy the scaling condition $\Nv(\constA\,\uh+\constB)=\constA^4\,\Nv$---{to achieve a complete dimensionless MSE};
		\item {it must be chosen such that the sampling error of the scale-invariant error estimator remains robust and fully invariant with respect to variations in the properties of the distribution of solution $\uh$.}
\end{enumerate}}


\rvsn{\textbf{$\var^2 \equiv {\mu_2}^2$ as a normalizer}: If one considers $\Nv = {\mu_2}^2$ as a standardizing factor, the MSE in \feq{Eq:MCerrorvar1} transforms to
\begin{equation}\label{Eq:MCerrorvarNMSE}
	\textrm{NMSE}(\htwoMC) = \frac{\Vtwo(\uh)}{{\mu_2}^2 N} +  \frac{({\var}(\uh)-{\var}(\uu))^2}{{\mu_2}^2}.
\end{equation}	
Provided that $0 < \Nv < \infty$, the normalization by $\mu_2^2$ enables the NMSE to be invariant under any linear transformation (as established in \feq{Eq:mu2square}). This satisfies the second requirement in the above list.
However, regarding the third criterion---robustness of the normalized sampling error to changes in the distribution of $\uh$---further analysis is required. The normalized sampling error, denoted by $\epsdot_s$, derived from \feq{Eq:V2}, is given as:
\begin{equation}\label{}
	\epsdot_s = \frac{1}{N} \left( \frac{\mu_4}{{\mu_2}^2} - \frac{N-3}{N-1} \right),
\end{equation} 
where we have suppressed the $\uh$ notation for clarity. In the asymptotic regime where $N \gg 1$, the ratio $(N - 3)/(N - 1) \to 1$, and the expression simplifies to:
\begin{equation}\label{eq:samperrnorm}
	\epsdot_s = \frac{1}{N} \left( \frac{\mu_4}{{\mu_2}^2} - 1 \right),
\end{equation}
in which the term $\mu_4 / \mu_2^2$ represents the standardized fourth central moment, kurtosis, denoted by $\kurt = {\mu_4}/{{\mu_2}^2}$ \cite{decarlo1997kurtosis}.
The normalized sampling error $\epsdot_s$ depends not only on the number of samples $N$ but also on the kurtosis of the solution distribution. This dependence becomes particularly significant for heavy-tailed distributions, which exhibit high kurtosis values. 
} 

\rvsn{
For example, in the case of a (symmetric) Gaussian distribution, the kurtosis is $\kurt = 3$ \cite{decarlo1997kurtosis}. In contrast, for a log-normal distribution, the kurtosis is given by $\kurt = \rho^4 + 2\rho^3 + 3\rho^2 - 3$ \cite{crow1988lognormal}, where $\rho = \exp(\var)$ and $\var$ is the variance of the underlying Gaussian distribution. If $\var = 1$, then $\kurt \approx 114$.
This implies that, for a fixed total NMSE and spatial resolution $\elesize$, the computational effort required to estimate the MC variance $\htwoMC$ is twice that of the mean estimator $\muMC$ for a Gaussian-distributed $\uh$ (since $\epsdot_s = 2/N$). However, for a log-normal distribution, the cost increases dramatically—by a factor of approximately 113 (i.e., $\epsdot_s \approx 113/N$)—relative to the mean estimator.
This violates the third criteria in the list, and therefore, undermines the effectiveness of $\var^2 \equiv {\mu_2}^2$ as a normalizer in practical applications.
}

Following this, we consider $\Nv = \Vtwo(\uh)$ as a normalizer for MC variance estimation. Given that $0<{\Vtwo}(\uh)<\infty$, and based on the scaling relation in \feq{Eq:V2scale}, it is evident that $\Vtwo(\uh)$ satisfies the second outlined requirement. 	
By normalizing the MSE with $\Vtwo(\uh)$ in \feq{Eq:MCerrorvar1}, we obtain the scale-invariant form of the MSE for the MC variance estimator:	
\begin{equation}\label{Eq:MCerror2b}
	\textrm{NMSE}(\htwoMC) :=\mathbb{E}({\e}^2) = \frac{1}{N} + \frac{ ({\var}(\uh)-{\var}(\uu))^2}{{\Vtwo(\uh)}}.
\end{equation}	
Here, the first term represents the scale-invariant sampling error, which decays proportionally to $1/N$ and is independent of the statistical characteristics of the solution $\uh$, thereby fulfilling the third requirement. The second term denotes the normalized squared discretization error.	

Furthermore, based on the computational complexity analysis for MC variance estimation discussed in \cite{krumscheid_quantifying_2020}, we henceforth recast the convergence analysis in terms of the normalized MSE formulation. 
\begin{proposition}
	Let us consider $\calpha, \cgamma, \alpha$ and $\convgamma$ as positive constants, and then one may define the error bounds, such that
	\begin{equation*}\label{eq:errorvarmc}
		\begin{split}
			(i)& \ \text{the deterministic error is bounded by}\ \frac{|{\var}(\uh)-\var(\uu)| }{\sqrt{\Nv}} \le \calpha \h^{\alpha},\\
			(ii)& \ \text{the cost to compute each sample of $\uh(x, \event)$ is given as} \\ 
			    & \mathcal{C}(\uh)  \le \cgamma \h^{-\convgamma}.
		\end{split}
	\end{equation*}
	Then, for any $0<\mathring{\epsilon}<e^{-1}$, the Monte Carlo (MC) variance estimator $\htwoMC$ with $\samplesize = \mathcal{O}(\mathring{\epsilon}^{-2})$ and $\elesize = \mathcal{O}(\mathring{\epsilon}^{1/\alpha})$ satisfies the normalized mean square error $\mathrm{NMSE}(\htwoMC) < \epssq$. Therefore, the corresponding computational cost of MC variance estimation is given as
	\begin{equation*}
		\mathcal{C}(\htwoMC)\le  \ccost \, \mathring{\epsilon}^{-2-\convgamma/\alpha},
	\end{equation*}
	where $\ccost>0$ is a positive constant. 
\end{proposition}

\rvsn{{\textbf{Unbiased estimation of $\Vtwo$}}: The quantity $\Vtwo(\uh)$ in \feq{Eq:MCerror2b} is analytical and not known. To determine the unbiased estimate of $\Vtwo(\uh)$, one must obtain unbiased estimates of entities $\mu_4(\uh)$ and $\mu_2(\uh)^2$, which is detailed in \app{appendix:unbiasMCvar}. Therefore, by substituting $\mu_4 \approx \text{h}_4^{\text{MC}}$ and $\mu_2^2 \approx \text{h}_{\left\lbrace 2,2\right\rbrace}^{\text{MC}}$ from \feqs{Eq:mu_4}{Eq:mu_2_2}, respectively, in \feq{Eq:V2}, one obtains the unbiased estimate: $\Nvhat:=\VtwoMC(\uh)$---of $\Vtwo(\uh)$. \rvst{A further analysis of the stochastic convergence of the estimator $\VtwoMC(\uh)$ is done in \app{appendix:varV2}.}
Finally, the MSE in \feq{Eq:MCerror2b} is re-described as, 
\begin{equation}\label{Eq:MCerror2c}
		\textrm{NMSE}(\htwoMC) :=\mathbb{E}({\e}^2) = \frac{1}{N} + \frac{ ({\var}(\uh)-{\var}(\uu))^2}{{\VtwoMC(\uh)}}.
\end{equation}
}

\section{Multilevel Monte Carlo method}\label{Sec:MLMC}

{The MSE's of MC mean and variance estimators in} \feqs{Eq:MCerror1}{Eq:MCerrorvar1}, respectively,
signify that to attain an overall higher level of total accuracy, one requires a very fine resolution of the \rvsb{finite element mesh} and a very large number of \acrshort{mc} samples. This demands a tremendous amount of computational effort, making the algorithm practically infeasible. Therefore, the desired moments are further estimated by a variance reduction technique in a multilevel fashion following \cite{Heinrich2001,giles2008,giles_2015}. 

\subsection{Multilevel Monte Carlo estimation of mean}\label{MLMC}

Let $\left\lbrace \meshlevel=0,1,2...,\LL \right\rbrace$ be a generalised increasing sequence---in the context of decreasing element size $h$---of nested \rvsb{meshes} $\mathcal{P}_\meshlevel$, a regular (non-degenerate) partition of the computational domain $\mathcal{G}$ of the problem described in \BRrefeq{abs_main_eq}. 
Here, $\meshlevel$ denotes the mesh level, and $\LL$ represents the finest mesh. 
The goal of the \acrshort{mlmc} method is to determine the statistics (such as the mean in this case) of the solution $\uhL(\vek{x}, \event)$ on the finest level $\LL$.
To this end, by exploiting the linearity of the expectation operator, one may express the \acrshort{mlmc} (for brevity, denoted by ML) mean estimate of the mean $\mumu(\uhL)$ using a set of samples $\left\lbrace \Nl \right\rbrace:=\left\lbrace \N_0,\N_1,\!...,\N_L \right\rbrace$ as \cite{giles_2015}
\begin{eqnarray}\label{Eq:meanMLMC}
	\muML(\uhLNL) &:=& \muMC(\uhONO) + \sum_{\meshlevel=1}^{\numofmesh}  \muMC(\uhlNl - \uhllNl) \nonumber\\
	&=& \sum_{\meshlevel=0}^{\numofmesh} \muMC(\Yl).
\end{eqnarray}
Here, $\muMC(\uhONO)$ is the \acrshort{mc} estimator of mean $\mu(\uhO)$ on level $l=0$ using $N_0$ samples, and $\muMC(\uhlNl - \uhllNl)$ represents the approximation of mean $\mu(\uhl - \uhll)$ with $\meshlevel>0$ and $\Nl$ samples. Furthermore, 
for $l=0$, $Y_0 = \uhONO$; else, $\Yl:=\uhlNl -\uhllNl$. Note that the individual term $\Yl, \ \meshlevel\ge{0}$, is sampled independently, and when $l>0$, the quantities $\uhlNl$ and $\uhllNl$ in $\Yl$ are considered to be strongly correlated---meaning that $\uhlNl$ and $\uhllNl$ are sampled from the same random seed.

As the mean estimate on the finest \rvsb{level} $\muML\equiv\muML(\uhLNL)$ is obtained as the telescopic sum of the difference of \acrshort{mc} mean estimates on the coarser \rvsb{levels}, the MSE of $\muML$, corresponding to \BRrefeq{Eq:MCerror1}, takes the form (see \cite{Cliffe2011}):
\begin{equation}\label{Eq:MLMCerror1}
	\textrm{MSE}(\muML) = \sum_{\meshlevel=0}^{\numofmesh} \frac{\var(\Yl)}{\Nl}+(\mathbb{E}(\uhL)-\mathbb{E}(u))^2.
\end{equation}
The above error consists of two terms: the variance of the estimator $\muML$ on the left and the square of the 
spatial discretization error on the right. 

\rvstt{\textbf{Scale-invariant {error estimator for MLMC mean estimation}}:} Similar to the MC error estimate for the mean in \feq{Eq:MCerror1}, the MSE of the MLMC mean estimator in \feq{Eq:MLMCerror1} depicts an absolute error estimate. As a result, the convergence of the MLMC algorithm also strongly depends on the solution magnitude. To this, based on the normalizer $\Nm$ defined in \fsec{MC}, and particularly the NMSE in \feq{Eq:MCerror3}, we consider $\Nm = \var(\uhL)$ as a standardizing entity, suggesting a new scale-invariant MSE estimate: 
\begin{equation}\label{Eq:MLMCerrornew2}
	\textrm{NMSE}(\muML) = \frac{1}{\var(\uhL)}\left(\sum_{\meshlevel=0}^{\numofmesh} \frac{\var(\Yl)}{\Nl}\right) +\frac{(\mathbb{E}(\uhL)-\mathbb{E}(u))^2}{\var(\uhL)}.
\end{equation}	
However, the normalization term $\var(\uhL)$, which is the variance of the solution $\uh$ on finest mesh level $L$, is generally not known. Hence, under the assumption that for any value of $l$, $\var(\uhl)$ will be approximately constant \cite{giles_2015,Cliffe2011}, one may substitute $\var(\uhL)$ with $\var(\uhO)$, which is based on the coarsest \rvsb{mesh} at $\meshlevel = 0$. Furthermore, under the consideration that $\Nmhat:= \htwoMC(\uhONO)$ is the MC estimator of $\var(\uhO)$, and by defining the MC variance estimate of $\var(\Yl):=\var(\uhl-\uhll)$ via $\htwoMC(\YhlNl)$, \feq{Eq:MLMCerrornew2} therefore transforms to	
\begin{equation}\label{Eq:MLMCerrornew3}
	\textrm{NMSE}(\muML) = \frac{1}{\htwoMC(\uhONO)}\left( \sum_{\meshlevel=0}^{\numofmesh} \frac{\htwoMC(\YhlNl)}{\Nl}\right)  +\frac{(\mathbb{E}(\uhL)-\mathbb{E}(u))^2}{\htwoMC(\uhONO)}.
\end{equation}
The first entity here stands for the estimated scale-invariant multi-level sampling error, whereas the second term defines the dimensionless squared discretization error. Therefore, to attain an overall normalized mean square error $\epssq$, it is adequate {that both terms are less than $\epssqtwo$.}
Note that the equal splitting of error $\epssq$
is not a requirement and can also be set otherwise---for more information on this, see \cite{collier_continuation_2015, haji-ali_optimization_2016}.

If $\costlevel\equiv\mathcal{C}(\YhlNl)$ is the computational cost of determining a single \acrshort{mc} sample of $\YhlNl$, then the overall cost of the estimator $\muML$ is given as $\mathcal{C}(\muML):=\sum_{\meshlevel=0}^{\LL}\Nl\costlevel$.
Here, the optimum number of samples $\Nl$ on each level $l$ is evaluated by solving an optimisation problem such that the {normalized sampling error is less than $\epssqtwo$.} 
As a result, the cost function
\begin{equation}\label{Eq:costfn}
	\costfunction = \text{arg} \min\limits_{\optSam}  \sum_{\meshlevel=0}^{\LL} \left(  \optSam\costlevel +  \Lagrange \frac{\htwoMC(\YhlNl)}{\Nmhat \Nl}\right)
\end{equation} 	
is minimised, due to which the optimal samples $\Nl$ are calculated as
\begin{equation}\label{eq:\Nl}
	\optSam = \Lagrange\left({\frac{\htwoMC(\YhlNl)}{\Nmhat \costlevel}}\right)^{\frac{1}{2}}.
\end{equation}
Here, $\Lagrange$ is the Lagrange multiplier, determined by
\begin{equation}\label{Eq:tau}
	\Lagrange = { {\frac{2}{\epssq}}}\sum_{\meshlevel=0}^{\numofmesh}\left( \frac{{\htwoMC(\YhlNl)\costlevel}}{\Nmhat}\right) ^{\frac{1}{2}}.
\end{equation} 	
\textcolor{black}{On the other hand, the set of mesh levels $\left\lbrace \meshlevel=0,1,2...,\LL \right\rbrace$ maybe optimally chosen (in a geometric or non-geometric sequence) for a given deterministic error. Typically, for PDE-based applications, the choice of coarse mesh $\meshlevel=0$ depends on the regularity of the solution $\uu(x)$ \cite{Cliffe2011}. Following this, the finer mesh selection is based on an apriori mesh convergence study, where the finest mesh \rvsb{level} $\meshlevel = L$ can be fixed or can be adaptively selected \cite{giles_2015} during the MLMC mean computation.  For a more detailed discussion on the optimal selection of mesh hierarchies, we refer the readers to \cite{haji-ali_optimization_2016}.}

The arguments for determining the total computational cost of the normalized MLMC mean estimator follow the classical MLMC procedure \cite{giles2008,Heinrich2001}. However, the difference is that the cost is described with respect to NMSE instead of MSE.
\begin{proposition}\label{prop:mlmcmean}
Let us consider the positive constants $\calpha, \cbeta, \cgamma, \convorder, \convbeta, \convgamma$, given that $\convorder\ge\frac{1}{2}\min(\convbeta,\convgamma)$. Then, we consider the following error bounds, where 
\begin{equation}\label{eq:errormean1}
	\begin{split}
		(i)& \ \text{the deterministic error decays as}\ \frac{|\mathbb{E}(\uhl)-\mathbb{E}(\uu)| }{\sqrt{\Nm}} \le \calpha \h_\meshlevel^{\convorder},\\
		(ii)& \ \text{the decay of variance is bounded by}\ \frac{\htwoMC(\YhlNl)}{\Nm} \le \cbeta \h_\meshlevel^{\convbeta}, \\
		(iii)& \ \text{the computational cost to determine a single realization of $\YhlNl$ is given as} \\ 
			 &\mathcal{C}(\YhlNl)  \le \cgamma \h_\meshlevel^{-\convgamma}.
	\end{split}
\end{equation}
Then, there exists another positive constant $\ccost$, such that for any $0<\mathring{\epsilon}<e^{-1}$, the multilevel Monte Carlo (MLMC) mean estimator satisfies the normalized mean square error $\mathrm{NMSE}(\muML)<\epssq$. Finally, the total computational cost of MLMC mean estimation is bounded by
\begin{equation}\label{Eq:costcomplex}
	\mathcal{C}(\muML)\le 
	\begin{cases}
		\ccost\,\mathring{\epsilon}^{-2}, \quad &\convbeta>\convgamma,\\
		\ccost\,\mathring{\epsilon}^{-2}(\log\mathring{\epsilon})^2, \quad &\convbeta=\convgamma,\\
		\ccost\,\mathring{\epsilon}^{-2-(\convgamma-\convbeta)/\convorder}, \quad &\convbeta<\convgamma.
	\end{cases}
\end{equation}
\end{proposition}
Based on the values of $\convbeta$ and $\convgamma$, one may also understand the major cost contributor amongst the sequence of mesh \rvsb{levels}. If $\convbeta>\convgamma$, the maximum cost is controlled by the coarsest level, and if $\convbeta<\convgamma$, the finest level governs the dominant cost. Finally, when $\convbeta = \convgamma$, the cost at each level is roughly evenly distributed.


Furthermore, it is clear for the first relation in \feq{eq:errormean1} that as $l\rightarrow \infty$,  $$\frac{|\mathbb{E}(\uhl)-\mathbb{E}(\uu)|}{\sqrt{\Nm}} \rightarrow 0.$$
However, $\mathbb{E}(\uu)$ is analytical and not known. Therefore, the deterministic error is defined via the triangle inequality as \cite{Cliffe2011}
\begin{equation}\label{Eq:trianglemean}
	\begin{split}
		\frac{|\mathbb{E}(\uhl) - \mathbb{E}(\uhll)|}{\sqrt{\Nm}}&\le \calpha \h_\meshlevel^{\convorder}.
	\end{split}
\end{equation}

%
%

\subsection{Multilevel Monte Carlo estimation of variance}\label{MLMCvar}
To enhance the characterization of the probabilistic system response $\uhL(\vek{x}, \event)$, this paper also focuses on defining the \acrshort{mlmc} estimator of $\var(\uhL)$ using h-statistics, expressed in the form \cite{krumscheid_quantifying_2020}:
\begin{eqnarray}\label{Eq:varMLMC}
	\htwoML(\uhLNL) &:=& \htwoMC(\uhONO) + \sum_{\meshlevel=1}^{\LL}\left( \htwoMC(\uhlNl) - \htwoMC(\uhllNl) \right).
\end{eqnarray}
Here, $\htwoMC(\uhONO)$ is the MC estimator of $\var(\uhO)$ with $\N_0$ samples; $\htwoMC(\uhlNl)$ and $\htwoMC(\uhllNl)$ represent the MC estimation of $\var(\uhl)$ and $\var(\uhll)$ using $\Nl$ samples, respectively.
For simplification, we introduce:
\begin{equation}\label{Eq:Zl}
	\Zl= 
	\begin{cases}
		\htwoMC(\uhONO),\quad &\meshlevel=0,\\
		\htwoMC(\uhlNl) - \htwoMC(\uhllNl),\quad &\meshlevel>0.
	\end{cases}
\end{equation}
Note that, for $\meshlevel>0$, $\uhlNl$ and $\uhllNl$ in $\Zl$ are determined using the same random seed.
Thus, the expansion in \feq{Eq:varMLMC} is rewritten as 
\begin{equation}
	\htwoML	= \sum_{\meshlevel=0}^{\numofmesh} \Zl.
\end{equation}
Similar to MLMC mean estimation, MLMC variance estimator is also obtained as the telescopic sum of the difference of \acrshort{mc} variance estimates on the coarser \rvsb{levels}. Therefore, in correspondence to \feq{Eq:MCerrorvar}, the MSE of the multilevel estimator $\htwoML$ takes the form:
\begin{equation}\label{Eq:MLMCerrorvar}
	\textrm{MSE}(\htwoML) = \var\left(\htwoML \right)  +  ({\var}(\uhL)-{\var}(\uu))^2.
\end{equation}
Under further consideration that the quantity $\Zl$, for $\meshlevel\ge0$, is sampled independently,
one may express the first term in the above equation as
\begin{equation}\label{Eq:sampleerror}
	\var\left(\htwoML \right)  = \sum_{\meshlevel=0}^{\numofmesh} \var(\Zl),
\end{equation}
in which the variance $\var(\Zl)$ is further defined---similar to \feq{Eq:V2N}---by ${\Vltwo}/{\Nl}$.
Therefore, \feq{Eq:MLMCerrorvar} is reformulated as
\begin{equation}\label{Eq:MLMCerrorvar1}
	\textrm{MSE}(\htwoML) = \sum_{\meshlevel=0}^{\numofmesh} \frac{\Vltwo}{\Nl}  +  ({\var}(\uhL)-{\var}(\uu))^2.
\end{equation}
Analogous to the MSE of MC variance estimator in \feq{Eq:MCerrorvar1}, here, the MSE is also split into statistical error---which is of the order $\mathcal{O}(\Nl^{-1})$---and square of deterministic bias. Furthermore, the quantity $\Vltwo$ is analytical, and the corresponding unbiased MC estimation, denoted $\VltwoMC$, is detailed in \app{appendix:unbiasMLvar}.		 

\rvstt{\textbf{Scale-invariant {error estimator for MLMC variance estimation}}:} The MSE of estimator $\htwoML$ in \feq{Eq:MLMCerrorvar1} is scale-dependent. Therefore, the MSE of multilevel variance estimator is transformed to a scale-invariant version $\textrm{NMSE}(\htwoML)$ by considering $\Nv = \Vtwo(\uhL)$ as the normalization entity---as described in \fsec{MCvar}.  Since, $\Vtwo(\uhL)$ is unknown, we make the assumption that the entity $\Vtwo(\uhl)$ stays approximately close for all values of $\meshlevel$. That is, the value of $\Vtwo(\uhL)$ on finest level $L$ is replaced by $\Vtwo(\uhO)$ of the coarsest mesh \rvsb{level} $\meshlevel = 0$. 
Finally, with the unbiased MC estimation of the normalizer $\Nvhat:=\VtwoMC(\uhONO)$, detailed in \fsec{MCvar} and \app{appendix:unbiasMCvar}, the new normalized MSE is given in the form:
\begin{equation}\label{Eq:MLMCerrorvar3}
	\textrm{NMSE}(\htwoML) = \frac{1}{\VtwoMC} \left( \sum_{\meshlevel=0}^{\numofmesh} \frac{\VltwoMC}{\Nl} \right)  +  \frac{({\var}(\uhL)-{\var}(\uu))^2}{\VtwoMC}.
\end{equation}			
Here, the first term defines the scale-invariant sampling error and the second term is the normalized squared deterministic error. The accomplishment of a total NMSE $\epssq$---as mentioned in \fsec{MLMC}---is justified by ensuring that both errors
are less than $\epssqtwo$.

One may determine the samples $\Nl$ analogously to \feeqs{Eq:costfn}{Eq:tau}, such that \feq{Eq:tau} is substituted in \feq{eq:\Nl}:		
\begin{equation}\label{}
	\optSam = { {\frac{2}{\epssq}}}\sum_{\meshlevel=0}^{\numofmesh}\left( \frac{{\VltwoMC\costlevel}}{\VtwoMC}\right) ^{\frac{1}{2}}\left({\frac{\VltwoMC}{\VtwoMC \costlevel}}\right)^{\frac{1}{2}}.
\end{equation}
Note that, here $\costlevel$ is the computational cost of evaluating one MC sample of $\uhlNl$ and $\uhllNl$ in the difference term $\Zl$, which is equivalent to determining the cost of $\YhlNl$ in \feq{Eq:costfn}.

\textcolor{black}{Following the computational complexity of the scale-invariant MLMC mean in \prop{prop:mlmcmean}, the computing cost of the normalized MLMC variance estimate, similar to the conventional MLMC variance (see \cite{krumscheid_quantifying_2020,bierig_convergence_2015}), is detailed below.}

\begin{proposition}\label{prop:mlmcvar}
	Let us introduce the positive constants $\calpha,  \cbeta, \cgamma, \convorderv, \convbetav, \convgamma$, such that $\convorderv\ge\frac{1}{2}\min(\convbetav,\convgamma)$. Then, one may define the following error bounds, such that
	\begin{equation}\label{eq:errorvar1}
		\begin{split}
			(i)& \ \text{the deterministic error is bounded by}\ \frac{|\var({u}_{h_l})-\var(u)|}{\sqrt{\Nv}} \le \calpha h_l^{\convorderv},\\
			(ii)& \ \text{the variance $\VltwoMC$ decays as}\ \frac{\VltwoMC}{\Nv}  \le \cbeta h_l^{\convbetav}, \\
			(iii)& \ \text{the computational cost to determine a single realization of $\Zl$ is given as} \\ 
			&\mathcal{C}(\Zl)  \le \cgamma \h_\meshlevel^{-\convgamma}.
		\end{split}
	\end{equation} 	
	Thus, one may state that, for any $0<\mathring{\epsilon}<e^{-1}$, the total normalized mean square error of multilevel Monte Carlo (MLMC) variance estimator  is bound by $\mathrm{NMSE}(\htwoML)<\epssq$. Finally, the computational cost of MLMC variance estimation reads
	\begin{equation}\label{Eq:costcomplex2}
		\mathcal{C}(\htwoML)\le 
		\begin{cases}
			\ccost\,\mathring{\epsilon}^{-2}, \quad &\convbetav>\convgamma,\\
			\ccost\,\mathring{\epsilon}^{-2}(\log\mathring{\epsilon})^2, \quad &\convbetav=\convgamma,\\
			\ccost\,\mathring{\epsilon}^{-2-(\convgamma-\convbetav)/\convorderv}, \quad &\convbetav<\convgamma.
		\end{cases}
	\end{equation}	
	Here, $\ccost>0$ is a positive constant. 
\end{proposition}

\textcolor{black}{Notably, if the deterministic and stochastic convergences in \feqs{eq:errormean1}{eq:errorvar1}} are correspondingly equal, then the computational complexity of the MLMC variance estimator will be asymptotically equal to that of the mean estimate. \textcolor{black}{However, this depends on the regularity of solution $\uu(x, \event)$, as proved by authors in \cite{bierig_convergence_2015}. Otherwise, generally, for a given total NMSE $\epssq$, the multilevel mean estimator runs faster as compared to the variance (elaborated with a numerical example in \fsec{NumericalResults}).}

Furthermore, as noted in \fsec{MLMC}, mesh level $\meshlevel=0$ contributes the most cost in the first scenario in the \feq{Eq:costcomplex2}; the second case imposes an almost equal cost on all mesh \rvsb{levels}, and the final relation demonstrates that the finest level $\LL$ is the most dominant.

Also, corresponding to \feq{eq:errorvar1}, when $l\rightarrow\infty$, there is a monotonic decay in ${|\var(\uhl)-\var(\uu)|}/{\sqrt{\Nv}}.$ Thereby, analogous to \feq{Eq:trianglemean}, the following condition
\begin{equation}\label{Eq:vardettri}
	\frac{|\var(\uhl)-\var(\uhll)|}{\sqrt{\Nv}} \le \calpha h_l^{\convorderv},
\end{equation}
holds.

\section{Application to linear elasticity: a model problem}\label{Sec:application}
\rvstt{To test the scale-invariant MC and MLMC methods, we consider a framework of linear elliptic PDEs with random coefficients. As a representative example, we focus on linear elasticity as a model problem, where the elasticity tensor is modelled as a matrix-valued random field capturing both the heterogeneous and random anisotropic nature of the material.}
\subsection{Deterministic setting}

Let $\spatialdomain\subset\Euclideanspace$ be a $d$-dimensional geometry with smooth Lipschitz boundary $\Gamma$. The aim is to determine the displacement vector $\dispvector \in \Euclideanspace$ (which belongs to the Hilbert space $\solspace$) at a spatial point $\vek{x}\in \spatialdomain$ that completely satisfy the equilibrium equations \cite{marsden_mathematical_1983}
\begin{align}\label{Eq:EquilForm}
	-\textrm{div} \, \stresstensor &= \volforcevector,  \quad \forall \vek{x} \in \spatialdomain, \nonumber \\
	\dispvector(\vek{x}) &= \vek{u}_{0}=\vek{0},  \quad \forall \vek{x} \in \Gamma_D, \\
	\stresstensor\cdot\vek{n}(\vek{x}) &= \tracvector, \quad \forall \vek{x} \in \Gamma_N, \nonumber
\end{align}
describing the linear-elastic behaviour. Here, $\stresstensor$ is the Cauchy stress tensor, which belongs to the space of second-order symmetric tensors $\SymSecRankSpace:=\left\lbrace \vek{\sigma}\in(\Euclideanspace\otimes\Euclideanspace) \ |\ \vek{\sigma=\vek{\sigma^T}}\right\rbrace $; $\volforcevector \in \Euclideanspace$ is the body force;  $\tracvector \in \Euclideanspace$ represents the surface tension on the Neumann boundary $\Gamma_N \subset \Gamma$, and $\vek{n}(\vek{x}) \in \Euclideanspace$ is the outward unit normal to $\Gamma_N$. For simplicity, a homogeneous boundary condition $\vek{u}_{0}=\vek{0} \in \Euclideanspace$ is considered on the Dirichlet boundary $\Gamma_D \subset \Gamma $. It is also possible to assume that $\Gamma_D \cap \Gamma_N = 0$.

One may further describe the strain-displacement relation as
\begin{equation}\label{Eq:kinematics}
	\straintensor = \frac{1}{2} \left( \nabla \dispvector + \nabla\dispvector^T \right), \quad  \forall \vek{x} \in \spatialdomain,
\end{equation}
where $\straintensor\in\SymSecRankSpace$ denotes an infinitesimal second-order symmetric strain tensor. Finally, the material constitutive equation is of the linear form
\begin{equation}\label{Eq:ConstLaw}
	\stresstensor = \elasttensor(\vek{x}) : \straintensor, \quad \forall \vek{x} \in \spatialdomain,
\end{equation}
in which $\elasttensor(\vek{x})$ represents a spatially varying fourth-order positive-definite symmetric elasticity tensor. Here, the notion of symmetry signifies the  major ($C_{{ijkl}} = C_{klij}$) and minor symmetries ($C_{ijkl} = C_{jikl} = C_{ijlk}$). As a result, one can map a $\elasttensor(\vek{x})$ tensor to a second-order tensor $\elastmatrix(\vek{x})$, more generally known as the elasticity matrix.  
The reduced $\elastmatrix(\vek{x})$ matrix belongs to a family of $(n\times n)$, where $n=d(d+1)/2$, real-valued positive-definite symmetric matrices:
\begin{equation}
\PosSymSecRankSpaceElast=\{ \elastmatrix\in(\mathbb{R}^n\otimes\mathbb{R}^n) \mid \elastmatrix=\elastmatrix^T, \vek{z}^T\elastmatrix \vek{z}>\vek{0}, \forall \vek{z}\in\mathbb{R}^n\setminus\vek{0}\}.
\end{equation}
Accordingly, conforming to Voigt notation, \feq{Eq:ConstLaw} transforms to $\sigma(\vek{x})=\elastmatrix(\vek{x})\cdot\varepsilon(\vek{x})$, with $\sigma(\vek{x})\in\Euclideanspace$ and $\varepsilon(\vek{x})\in\Euclideanspace$ denoting the stress and strain vectors, respectively.

\subsection{{Stochastic modelling of material uncertainty: a reduced parametric approach}}\label{sec:stomodel}

Material law as described previously is complicated when highly heterogeneous and anisotropic materials---such as bone tissue, see \fsec{NumericalResults}---are to be modelled. To include aleatoric uncertainty with heterogeneity in \feq{Eq:ConstLaw}, the material properties have to be modelled as random and spatially dependent \cite{sobczyk_stochastic_2001, ostoja-starzewski_microstructural_2007}. In this paper, a probabilistic point of view is studied, in which the $\elastmatrix(\vek{x})$ matrix is modelled as a matrix-valued second-order random field $\CRF$ on a probability space $(\samplespace,\mathcal{F},\mathbb{P})$. In other words, the random elasticity matrix can be modelled as a mapping:
\begin{equation}
	\CRF:\spatialdomain\times\samplespace\rightarrow\PosSymSecRankSpaceElast. 
\end{equation}  
Furthermore, in practise, the usual assumption is that the elasticity matrix $\elastmatrix(\vek{x})$ globally follows a certain type of symmetry (e.g., isotropic, orthotropic, etc.; see \cite{cowin_continuum_2013,bona_coordinate-free_2007} for invariance classification of elasticity matrix $\elastmatrix$), even though the experiments today cannot provide this information with full certainty. 
Therefore, to include uncertainty in the type of material symmetry, in this paper, we model the elasticity matrix $\elastmatrix(\vek{x})$  as per the reduced parametric approach (also commonly known as the non-parametric method) \cite{soize2006non,soize2008,GUILLEMINOT2010}. 
That is, the random matrix $\CRF$ follows a specific type of symmetry only in the mean, whereas each of the realisations belongs to the triclinic system, which is the lowest order of symmetry for elasticity-type tensors. 
This gives our model a full degree of freedom in a case where the predefined mean symmetry turns out to be incorrect; one would still be able to identify other types of invariances given experimental data. This is, however, not the case if we assume that each of the realisations is constrained. \rvsn{If the specific symmetry class of the material---beyond the triclinic case---is known for both the ensemble (population) and its mean, one may employ the stochastic modelling frameworks presented in \cite{malyarenko_tensor-valued_2019}, and \cite{shivanand_stochastic_2024, shivanand_stochastic_mms}, which also has the ability to separate the modelling of strength, eigen-strain distribution, and spatial orientation, allowing for independent control of each component; however, this remains beyond the scope of the present study.}
 

Following this, we model the homogeneous mean matrix as $\mathbb{E}(\CRF)=\uppermean^T\uppermean$,
in which the term $\uppermean \in \mathbb{R}^{n\times n}$ represents an upper triangular matrix---this square-type (Cholesky) factorization ensures the positive-definiteness of the mean matrix. 
Then, to allow uncertainty into the model, the mean formulation is extended to 
\begin{equation}\label{Eq:RandomCmatrix}
	\CRF = \uppermean^T\stogermnoX\uppermean,
\end{equation}
such that
\begin{equation}\label{Eq:IdentityT}
	\mathbb{E}({\stogermnoX}) = \identitymatrix
\end{equation}
holds. Here, $\stogermnoX$ is a matrix-valued random variable that also resides in $\PosSymSecRankSpaceElast$, the mean value of which is an identity matrix $\identitymatrix\in\PosSymSecRankSpaceElast$. In this manner, the mean behaviour is controlled by $\mathbb{E}(\CRF)$, and the random fluctuations are governed by $\stogermnoX$. To construct a random ensemble $\stogermnoX$, one may use the maximum entropy optimisation principle under the constraints mentioned in the previous list (but with a mean defined in \feq{Eq:IdentityT}). 
Henceforth, $\stogermnoX$ can further be factorised as
\begin{equation}\label{Eq:choleskyT}
	\stogermnoX = \upperT^T\upperT,
\end{equation}
and by substitution in \feq{Eq:RandomCmatrix}, one obtains
\begin{equation}\label{Eq:Randommatrix}
	\CRF = \uppermean^T\upperT^T\upperT\uppermean.
\end{equation}
Here, $\upperT\in\mathbb{R}^{n\times n}$ denotes the upper triangular random matrix with entries
\begin{equation}\label{Eq:upperelements}
	{V}_{ij}= 
	\begin{cases}
		\frac{\coeffdispT}{\sqrt{d+1}}\theta_{ij}(\vek{x},\event), \text{if} \  i<j\\
		\frac{\coeffdispT}{\sqrt{d+1}}\sqrt{\normconst\GammaRF}, \text{if} \ i=j,
	\end{cases}
\end{equation}
given that the elements ${V}_{ij}, i\le j$ are independent.
Observe that the non-diagonal upper triangular elements are modelled as 
independent Gaussian random fields $\theta_{ij}(\vek{x},\event), 1\leq i\leq j \leq n$ with zero mean and unit variance. On the other hand, the diagonal entries are gamma-distributed positive-definite random fields:
\begin{gather}\label{Eq:inverseCDF}
	\GammaRF = \gammaCDF^{-1}\circ \erf(\theta_{ij}(\vek{x},\event)).
\end{gather}
Here,
\begin{equation}
	\alpha_j = (d+1)/(2\coeffdispT^2)+(1-j)/2
\end{equation}
is a positive real number; $\erf$ is the standard Gaussian distribution function; and $\gammaCDF^{-1}$ is the inverse gamma cumulative distribution function. This assures that the diagonal elements are strictly positive, and therefore, the random matrix $\stogermnoX$ also remains positive-definite.
Reverting to the above description in \feq{Eq:upperelements}, where $\normconst$ is a normalization constant of
the resultant function (given in \feq{Eq:inverseCDF}), the value of which equals 2.
Furthermore, $\coeffdispT:=[0,1]\in\mathbb{R}$ defined as \feq{Eq:delT}, is a scalar value that controls the dispersion of $\stogermnoX$:
\begin{equation}\label{Eq:delT}
	\coeffdispT=\left\lbrace \frac{1}{d}\mathbb{E}\left[ \parallel \stogermnoX - {\identitymatrix}\parallel^2 \right] \right\rbrace^{1/2}.
\end{equation}
The coefficient of dispersion parameter $\coeffdispT$ is chosen such that 
\begin{equation}\label{Eq:deltaC}
	\coeffdispC = \left[ \frac{\mathbb{E}\left\lbrace \lVert \CRF - \mathbb{E}(\CRF) \rVert^2 \right\rbrace}{\lVert \mathbb{E}(\CRF) \rVert^2}\right]^{1/2} = \frac{\coeffdispT}{ \sqrt{d+1}} \left[ 1+\frac{(\text{tr}(\mathbb{E}(\CRF)))^2}{\text{tr}(\mathbb{E}(\CRF))^2} \right]^{1/2}
\end{equation}
holds. Here, $\coeffdispC:=[0,1]\in\mathbb{R}$ is the coefficient of dispersion of the random matrix $\CRF$.

\subsection{Stochastic setting}\label{SubSec:stosetelast}

The description of the random elasticity matrix field $\CRF$ leads to the apparent transformation of a linear-elastic material model to a stochastic model. The aim is to determine the random displacement vector field $\dispvector(\vek{x},\event):\spatialdomain\times\samplespace\rightarrow\Euclideanspace$. Therefore, the equilibrium equations are rewritten in the form:
\begin{align}\label{Eq:stoEquilForm}
	-\textrm{div} \, \vek{\sigma}(\vek{x},\event) &= \volforcevector,  \quad \forall \vek{x} \in \spatialdomain, \event\in\samplespace, \nonumber \\
	\dispvector(\vek{x},\event) &= \vek{u}_{0}=\vek{0},  \quad \forall \vek{x} \in \Gamma_D, \event\in\samplespace, \\
	\vek{\sigma}(\vek{x},\event)\cdot\vek{n}(\vek{x}) &= \tracvector, \quad \forall \vek{x} \in \Gamma_N, \event\in\samplespace, \nonumber
\end{align}
in which the boundary conditions and body forces $\vek{f}(\vek{x})$ remain deterministic. 
Further, the linearized kinematics relationship is transformed to
\begin{equation}
	\vek{\varepsilon}(\vek{x},\event) = \frac{1}{2} \left( \nabla \dispvector(\vek{x},\event) + \nabla\dispvector(\vek{x},\event)^T \right), \forall \vek{x} \in \spatialdomain, \event\in\samplespace,
\end{equation}
where the $\nabla(\cdot)$ operator is taken in a weak sense. Finally, 
the constitutive law is represented by
\begin{equation}\label{Eq:stoConstLaw}
	\vek{\sigma}(\vek{x},\event) = \elasttensor(\vek{x},\event) : \vek{\varepsilon}(\vek{x},\event), \forall \vek{x} \in \spatialdomain, \event\in\samplespace.
\end{equation}
By carrying out the variational formulation of the above stochastic partial differential equations on $\spatialdomain$ and further discretizing in a finite element setting \cite{zienkiewicz_finite_2013}, one searches for the solution $\dispvectorApprox(\vek{x},\event):\finitespace\times\samplespace\rightarrow\Euclideanspace$ in a finite subspace, $\finitespace$, as described in \cite{lord_powell_shardlow_2014, rosic_variational_2012,Keese_thesis}. 

\rvsb{In computational stochastic mechanics, a vast body of literature exists on numerical methods for obtaining stochastic solution $\dispvectorApprox(\vek{x},\event)$ after semi-discretization. Some well-known approaches that fall into the category of series-expansion methods include the spectral stochastic finite element method \cite{ghanem_stochastic_1991,matthies_uncertainties_1997,matthies_finite_1999}, the perturbation method \cite{liu_probabilistic_1986,Kleiber1992,kaminski_generalized_2007}, and the Neumann expansion method \cite{yamazaki_neumann_1988,shinozuka_response_1988}. 
The present study, however, focuses on another class of techniques that directly integrate the statistics of the response, of which, Monte Carlo simulation \cite{Shinozuka1995,papadrakakis_robust_1996} has traditionally been used for assessing the validity of other methods. 
Further reviews and discussions of the procedures can be found in \cite{ghanem_stochastic_1991, matthies_uncertainties_1997, schueller_computational_2001,stefanou_stochastic_2009}. 
}

The goal of this study is to determine the total displacement scalar field in Euclidean norm, i.e., $\totaldisp(\vek{x}, \event)=\|\dispvectorApprox(\vek{x}, \event)\|$, and estimate the second-order statistics such as the mean and variance of sampled response $\uhNt:=[\totaldisp(\vek{x},\event_i)]_{i=1}^\N$ using the scale-invariant MC and MLMC as detailed in \fsecs{Sec:MC}{Sec:MLMC}. 
\rvsb{Note that a comparison of the proposed MC and MLMC techniques with the previously listed series-expansion approaches lies outside the scope of this work.}

\newlength\figH
\newlength\figW
\setlength{\figH}{3.9cm}
\setlength{\figW}{3.9cm}

\newlength\figHh
\newlength\figWw
\setlength{\figHh}{6.2cm}
\setlength{\figWw}{7.5cm}

\definecolor{mycolor1}{rgb}{0.93,0.69,0.13}%
\definecolor{mycolor2}{rgb}{0.64,0.08,0.18}%

\newlength\figHH
\newlength\figWW
\setlength{\figHH}{4cm}
\setlength{\figWW}{4cm}

\newlength\figHNT
\newlength\figWNT
\setlength{\figHNT}{3.2cm}
\setlength{\figWNT}{3.2cm}

\section{Numerical results: 2D human femur}\label{NumericalResults}

\rvstt{In this section, we examine the performance of the scale-invariant MC and MLMC algorithms on a linear elastic material model of a 2D human femoral bone, which is a highly heterogeneous and anisotropic material.}

\subsection{Specifications}
\begin{figure}[h]
	\centering
	\includegraphics[width=0.24\textwidth]{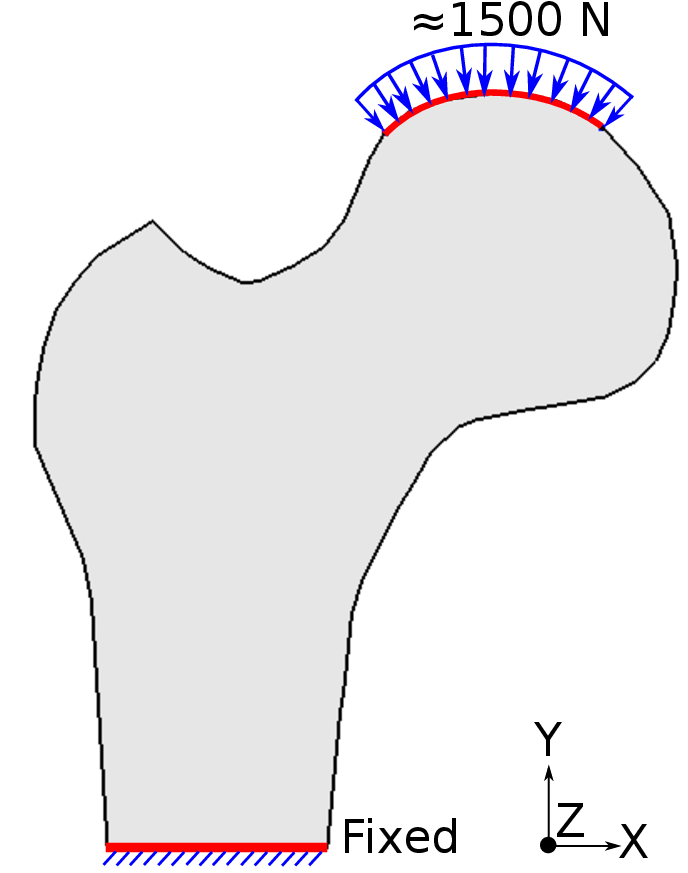}
	\caption{Geometry and boundary conditions \label{fig:femur_bc}}
\end{figure}

A two-dimensional proximal femur bone geometry with a body width of approximately \rvst{7 cm and 21.7 cm} in total height is considered.   \ffig{fig:femur_bc} shows the boundary conditions where an in-plane uniform pressure load with a resultant load of \rvst{1500 N} is applied on top of the bone and zero displacements are considered at the bottom \cite{yosibash_ct-based_2007}.
The \acrfull{fem} based spatial discretization is undergone using four-noded plane stress elements.
By sampling the probabilistic space, each deterministic simulation is thus executed by the finite element MATLAB-based software Plaston \cite{rosic_variational_2012}, \rvsb{where preconditioned conjugate gradient (PCG) method is used as an iterative solver.}

\begin{figure}[ht]
	\centering
	\includegraphics[width=1.0\textwidth]{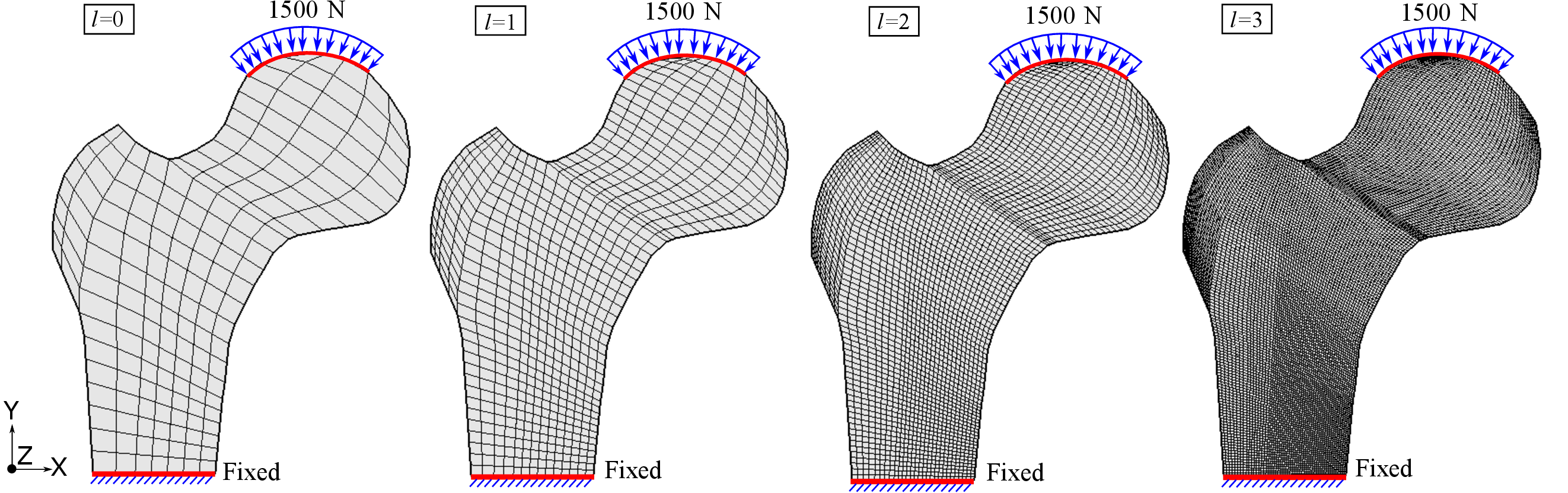}
	\caption{Nested mesh \rvsb{levels} of 2D femur bone}
	\label{Fig:4_levels}
\end{figure}

To implement the scale-invariant \acrshort{mlmc} method, a sequence of four nested meshes with element size $\h_{\meshlevel-1}=2\h_\meshlevel$ are considered, as shown in \ffig{Fig:4_levels}.
One may notice the identical implementation of boundary conditions over all the mesh levels.
The corresponding mesh specifications are tabulated in \ftbl{Table:mesh_spec}.
\begin{table}[ht]
	\centering 
	{\begin{tabular}{llll}
			\specialrule{1pt}{1pt}{1pt}			
			$l$
			& Elements & Nodes & DOF \\ \midrule
			0 & 171 & 206 & 396 \\ 
			1 & 684 & 753 & 1476 \\ 
			2 & 2736 & 2873 & 5688 \\  
			3 & 10944 & 11217 & 22320 \\ 
			\specialrule{1pt}{1pt}{1pt}
	\end{tabular}}
	\caption{Mesh specifications of 2D meshes\label{Table:mesh_spec}}
\end{table}

\rvsn{Femoral bone tissue is not only a highly heterogeneous material but also anisotropic, with its precise type of elastic symmetry remaining uncertain \cite{Sansalone:2016,KELLER19941159}. To account for this, we model the material stiffness using a matrix-valued random field $\CRF$, as defined in \fsec{sec:stomodel}.}
\begin{table}[ht]
	\centering
	{\begin{tabular}{lll}
			\specialrule{1pt}{1pt}{1pt}
			\begin{tabular}[c]{@{}l@{}}
				Young's \\ modulus (GPa) 
			\end{tabular} & 
			\begin{tabular}[c]{@{}l@{}} 
				Poisson's \\ ratio 
			\end{tabular} & 
			\begin{tabular}[c]{@{}l@{}} 
				Shear\\ modulus (GPa) 
			\end{tabular} \\ \midrule
			$E_1 = 12 $ & $\nu_{21} = 0.371$ & $G_{12} = 5.61$\\ 
			$E_2 = 20$ &  & \\
			\specialrule{1pt}{1pt}{1pt}
	\end{tabular}}
	\caption{Orthotropic material parameters\label{Table:mat_param_orth}}
\end{table}
\rvsn{We consider the average orthotropic elastic properties of human cortical femoral bone, as reported in the experimental study conducted on 60 specimens by \cite{ashman_continuous_1984}. The corresponding elastic coefficients, restricted to 2D, are listed in \ftbl{Table:mat_param_orth}, and are used to define the homogeneous mean matrix $\mathbb{E}(\CRF)$ of the random elasticity tensor field $\CRF$.} 
By setting the coefficient of dispersion of the random elasticity matrix field to $\coeffdispC = 0.1$, the coefficient of dispersion of the matrix-valued random field $\stogermnoX$ is determined using \feq{Eq:deltaC}.
\begin{wrapfigure}[13]{r}{0.39\textwidth}
	\centering
	\input{images/KL_eigenvals}
	\caption{\rvsn{Decay of eigenvalues across all mesh levels}}
	\label{fig:decayeigs}
\end{wrapfigure}
Accordingly, the fluctuation matrix $\stogermnoX$ is modelled as a non-linear transformation of 6 (as $n=3,\ d=2$) independent scalar Gaussian random fields. Each of these is approximated via a truncated (up to $\KLEterms$ terms) Kosambi-Karhunen-Lo\`{e}ve expansion \cite{Kosambi_1943,Karhunen1947,Loeve1948}, expressed as: 
\rvsn{\begin{equation}\label{Eq:GRF}
	\theta_{ij}(\vek{x},\event) = \overline{\theta}_{ij}(\vek{x}) + \sum_{k=1}^{\KLEterms} \sqrt{\eta_k}\psi_k(\vek{x}) \xi_k(\event).
\end{equation}
Here, $\overline{\theta}_{ij}(\vek{x})$ denotes the spatially varying mean field, $\xi_k(\event)$ are the mutually uncorrelated and independent standard Gaussian random variables, and $(\eta_k,\psi_k)$ are the eigenpairs of the covariance operator associated with the autocorrelation kernel  $\left\lbrace R_{ij}(r) : R_{ij}(0)=1,1\leq i \leq j \leq n\right\rbrace$. These eigenpairs are obtained via a finite element discretization of the Fredholm integral equation of the second kind \cite{spanos_stochastic_1989, Keese_thesis, matthies_galerkin_2005, betz_numerical_2014}. 
In this work, we assume a Gaussian-type covariance structure:  
\begin{equation}\label{Eq:corlength}
	R_{ij}(r) = \varrho^2 \exp(-\corrlength^{-2}r^2),
\end{equation}
in which, $r=\| \vek{x}_1-\vek{x}_2\|$ is the Euclidean distance between spatial points $(\vek{x}_1,\vek{x}_2)\in\spatialdomain$ with $r\geq0$, and $\varrho^2$ is the marginal variance. Moreover, each of $R_{ij}$ are individually parametrized by a vector of correlation lengths $l_c\in\Euclideanspace_+$. For simplicity, we adopt identical correlation lengths across all components of autocorrelation $R_{ij}$.} \rvsb{Note that the standard Gaussian random variables $\xi_k(\event)$ in \feq{Eq:GRF} are generated using MATLAB's built-in \textit{randn} function, which by default employs the ziggurat algorithm \cite{marsaglia_fast_1984,marsaglia_ziggurat_2000}.
} 

In this study, the mean $\overline{\theta}_{ij}(\vek{x})$ is set to be spatially constant with a value of zero and the variance $\varrho^2=1$.
\rvst{Under the assumptions considered for constructing the MLMC estimator in \fsec{MLMC}, each matrix-valued random field $\CRF$ and the corresponding Gaussian fields $\theta_{ij}(\vek{x},\event)$ are modelled independently on each mesh level $\meshlevel$. Whereas, for the definition of the difference terms $\Yl$ and $\Zl$ in \feqs{Eq:meanMLMC}{Eq:Zl}, the random field on the coarse mesh level $l-1$ is obtained by directly mapping it from the fine level $l$ at the intersecting common spatial nodes.}
Furthermore, the Gaussian autocorrelation function in \feq{Eq:corlength} is defined with a correlation length of \rvst{3.5 cm} in both $x$- and $y$-directions, across all four levels of meshes.
The expansion, given in \feq{Eq:GRF}, is truncated to $\KLEterms = 100$ terms \rvst{on all mesh levels}, chosen based on the decay of eigenvalues $\eta_k$ \rvst{on fine mesh $L$}, as illustrated in \ffig{fig:decayeigs}. \rvst{For further optimality, one may also use level-dependent truncations, as discussed in \cite{teckentrup_further_2013}.} 

\begin{figure}[h!]
	\centering
	\begin{minipage}{0.02\textwidth}
		\vspace{-3.5cm} 
		\rotatebox{90}{{\footnotesize Mean (Orthotropic)}}
	\end{minipage}%
	\begin{subfigure}{0.31\textwidth}
		\centering
		\includegraphics[width=0.75\textwidth]{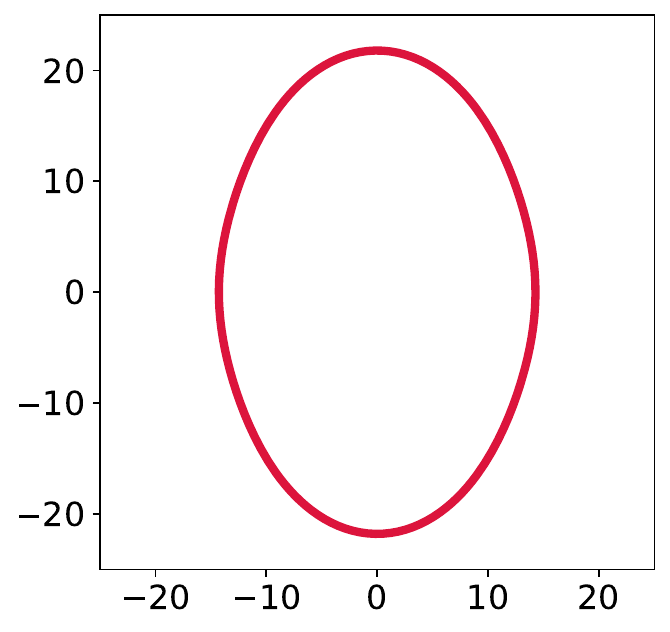}
	\end{subfigure} 
	\begin{subfigure}{0.31\textwidth}
		\centering
		\includegraphics[width=0.74\textwidth]{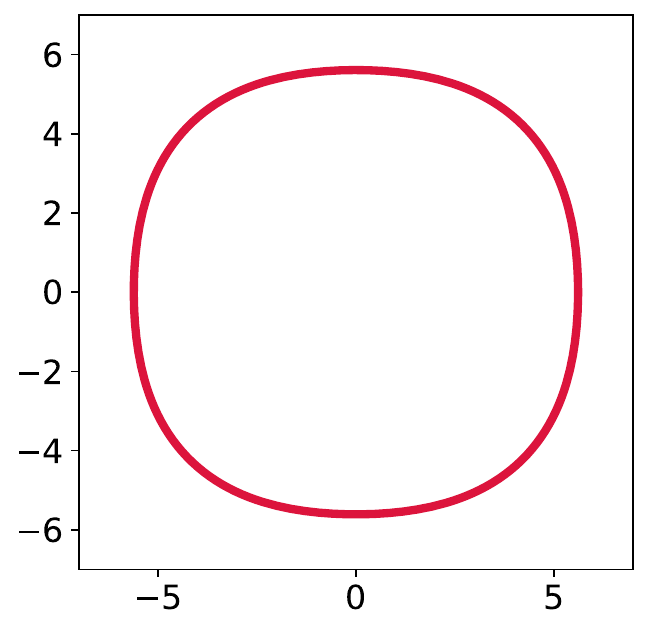}
	\end{subfigure}
	\begin{subfigure}{0.31\textwidth}
		\centering
		\includegraphics[width=0.75\textwidth]{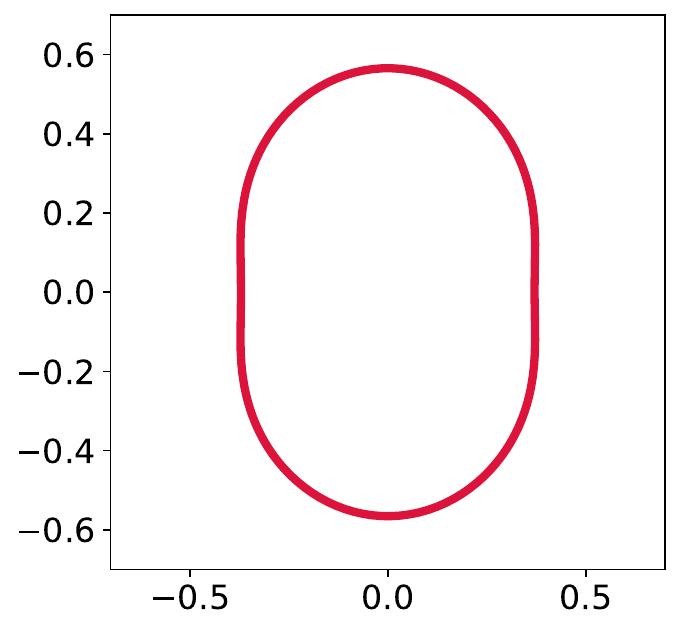}
	\end{subfigure}
	\begin{minipage}{0.02\textwidth}
		\vspace{-3.5cm} 
		\rotatebox{90}{{\footnotesize Sample 1 (Triclinic)}}
	\end{minipage}%
	\begin{subfigure}{0.31\textwidth}
		\centering
		\includegraphics[width=0.75\textwidth]{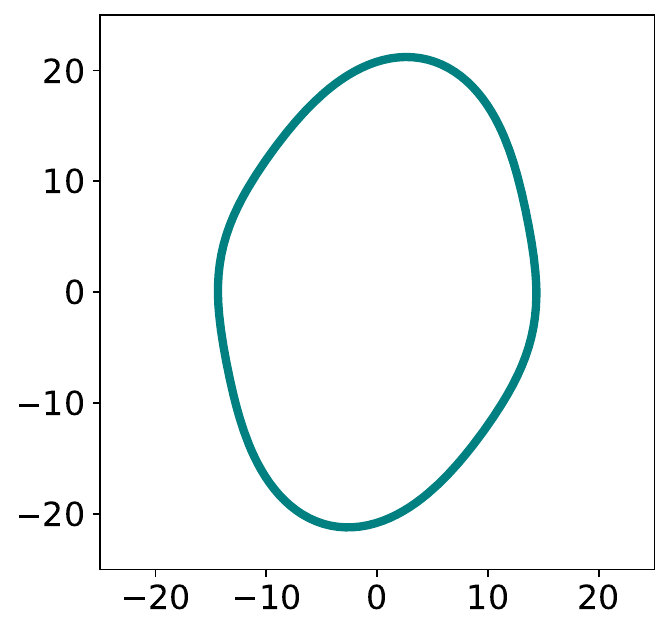}
	\end{subfigure} 
	\begin{subfigure}{0.31\textwidth}
		\centering
		\includegraphics[width=0.74\textwidth]{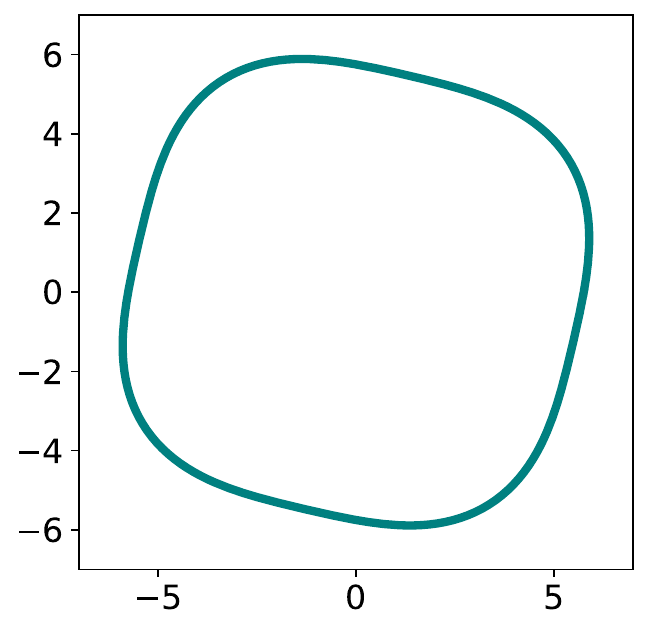}
	\end{subfigure}
	\begin{subfigure}{0.31\textwidth}
		\centering
		\includegraphics[width=0.75\textwidth]{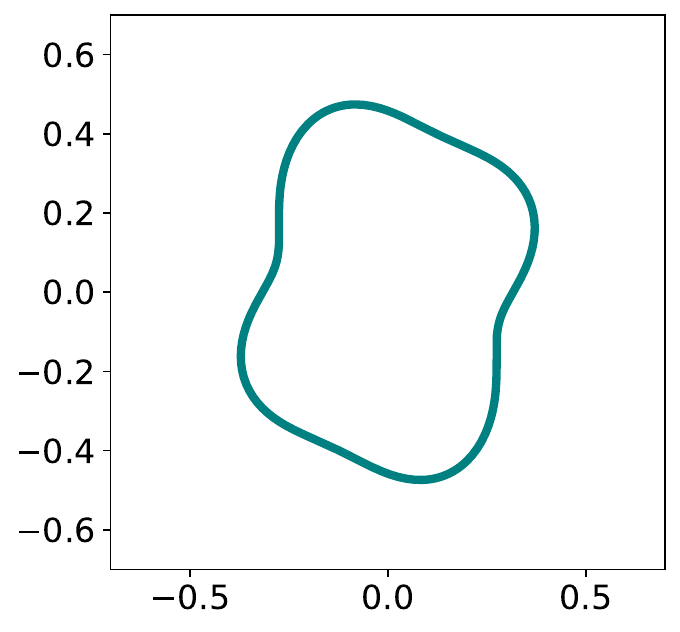}
	\end{subfigure}
	\begin{minipage}{0.02\textwidth}
		\vspace{-4.5cm} 
		\rotatebox{90}{{\footnotesize Sample 2 (Triclinic)}}
	\end{minipage}%
	\begin{subfigure}{0.31\textwidth}
		\centering
		\includegraphics[width=0.75\textwidth]{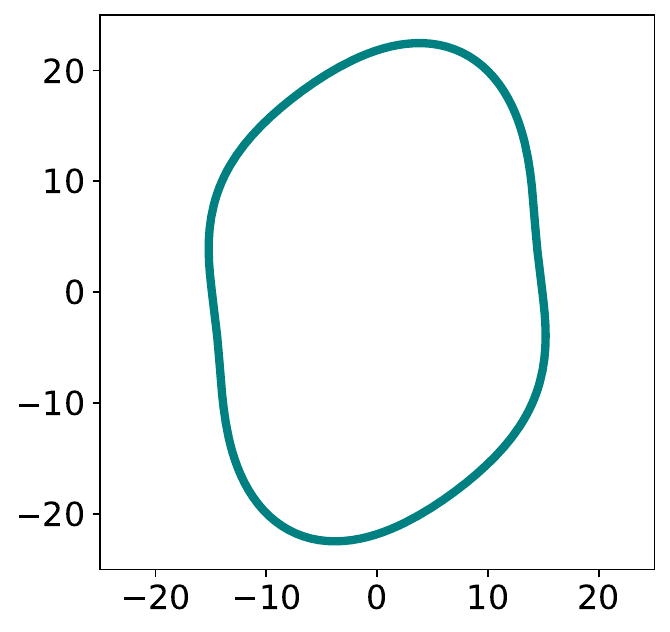}
		\caption{Young's modulus}
	\end{subfigure} 
	\begin{subfigure}{0.31\textwidth}
		\centering
		\includegraphics[width=0.74\textwidth]{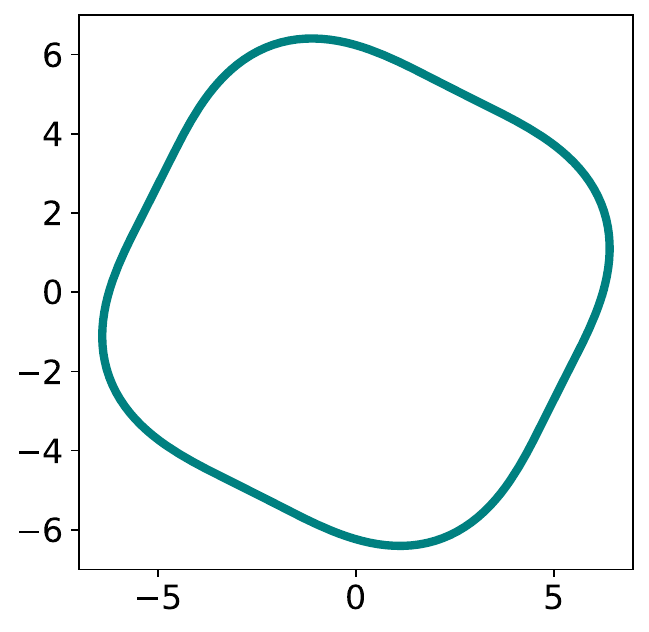}
		\caption{Shear modulus}
	\end{subfigure}
	\begin{subfigure}{0.31\textwidth}
		\centering
		\includegraphics[width=0.75\textwidth]{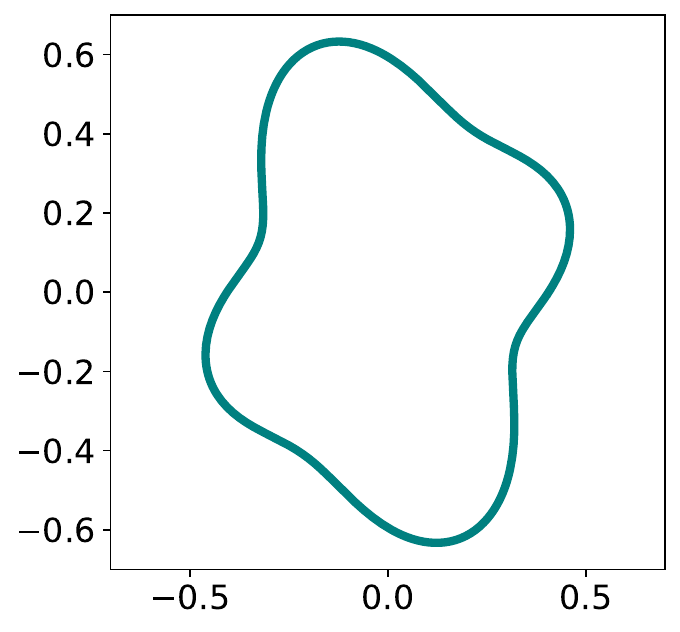}
		\caption{Poisson's ratio}
	\end{subfigure}
	\caption{\rvsn{Visualization of random anisotropy $C(\cdot, \event)$ on mesh $\meshlevel=0$}}
	\label{fig:randaniso}
\end{figure}

\begin{figure}[h!]
	\centering
	\includegraphics[width=1.0\textwidth]{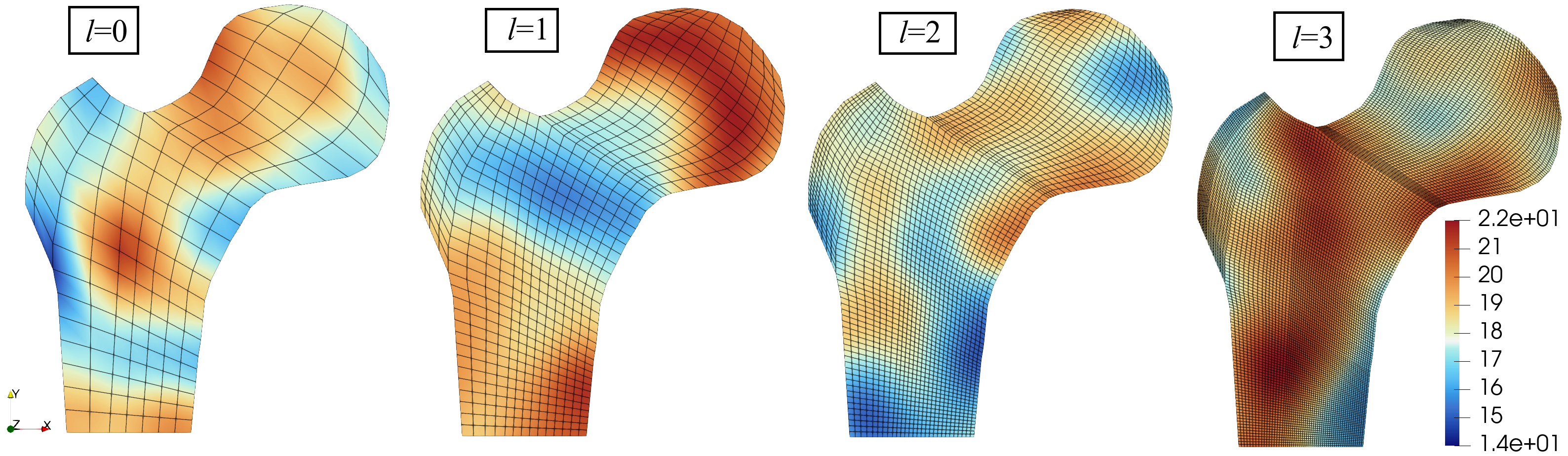}
	\caption{\rvsn{Visualization of spatial variation of a realization $C_{1, 1}(\vek{x}, \cdot)$}}
	\label{fig:spatvar}
\end{figure}

\rvsn{The random anisotropy of the material at a fixed spatial location $\vek{x}$ on the coarse mesh \rvsb{level} $\meshlevel=0$, modelled by the matrix-valued random field $C(\cdot, \event)$, is illustrated in \ffig{fig:randaniso}. Here we demonstrate the characteristic directional elastic parameters in the columns---namely, Young’s modulus, shear modulus, and Poisson’s ratio—using the open-source software ELATE \cite{gaillac_elate_2016}. The first row of the figure presents the orthotropic elastic properties corresponding to the mean matrix $\mathbb{E}(\CRF)$, while the subsequent rows display triclinic characteristics of two individual realizations, $C(\cdot, \event_1)$ and $C(\cdot, \event_2)$, where one may notice the variation in shape and size of all three parameters as compared to the mean. Additionally, we visualize the spatial variation of a single realization of the random field component $C_{1,1}(\vek{x}, \cdot)$ across all mesh levels in \ffig{fig:spatvar}.}


The objective of the study is to compute the \acrshort{mlmc} mean $\muML$ and variance $\htwoML$ estimates of the total displacement random field ${u}_{h_L}^{(t)}(\vek{x}, \event)$. The procedure for implementation of the scale-invariant MLMC is similar to that of the conventional MLMC method; the primary difference is in the usage of normalized error instead of absolute error; see \cite{giles_2015} for the algorithm. 
In this study, we assume that the optimal finest level (here, $L=3$) is known i.e., {the normalized squared discretization errors, from \feqs{Eq:MLMCerrornew3}{Eq:MLMCerrorvar3}, are less than $\epssqtwo$.} 
Therefore, the focus is only on satisfying the condition of normalized sampling errors; in other words, the objective is to ensure that
{the normalized sampling errors in \feqs{Eq:MLMCerrornew3}{Eq:MLMCerrorvar3} are less than $\epssqtwo$.}
Furthermore, to avoid the interpolation error, the terms $\YhlNl$ and $\Zl$ in \feqs{Eq:meanMLMC}{Eq:Zl} are calculated only at finite element nodes that have the same common spatial coordinates between all four levels of meshes. This further means that the MLMC mean and variance estimates of the system response on the finest level $L$ are evaluated only at these common nodes.

\subsection{Screening test}


An \rvst{a priori} performance analysis of MLMC mean and variance estimators is done by the so-called screening test, in which a fixed number of 50 samples is considered over four levels of meshes. 
For the MLMC mean estimate, we assume---for \rvst{simplicity}---that the estimated normalizing function $\htwoMC(\uhONOt)$ (from \feq{Eq:MLMCerrornew3}) is spatially constant, i.e., $\Nmhatinf =\max(\htwoMC(\uhONOt))$. With this, by considering a ${L}_\infty$ norm over the error bounds in \feq{Eq:trianglemean} and the second relation in \feq{eq:errormean1}, one obtains: 
\begin{equation}\label{Eq:convNT}
	\begin{split}
		\frac{\max|\muMC(\YhlNl)|}{\sqrt{\Nmhatinf}} &\le c_{\alpha} \h_\meshlevel^{\convorder},\\
		\frac{\max\left(\htwoMC(\YhlNl)\right)}{\Nmhatinf}&\le c_{\beta} \h_\meshlevel^{\convbeta}.
	\end{split}
\end{equation}
As to the MLMC variance estimate, similar to the assumptions made in the previous equations, one may consider $\Nvhatinf=\max(\VtwoMC(\uhONOt))$ as the normalization constant. Thereby, \feq{Eq:vardettri} and the second case in \feq{eq:errorvar1} are rewritten as
\begin{equation}\label{Eq:convNTvar}
	\begin{split}
		\frac{\max|\Zl|}{\sqrt{\Nvhatinf}} &\le c_{\alpha} h_l^{\convorderv},\\
		\frac{\max\left(\VltwoMC\right)}{\Nvhatinf} &\le c_{\beta} h_l^{\convbetav}.
	\end{split}
\end{equation} 	

Following this, \ffig{Fig:meanscreening} shows an overview of the corresponding results.
Note that the normalization constants $\Nmhatinf$ and $\Nvhatinf$ for the screening test are also set by 50 initial samples; \rvst{however, during the implementation of the MLMC method, the estimates of both constants are progressively updated as additional samples become available}.
The top left and right plots show the behaviour of the logarithm of ratios defined in  \feqs{Eq:convNT}{Eq:convNTvar}.
The deterministic decay of difference terms $\max|\muMC(\Yl)|$ and $\max(\Zl)$ can be seen. A stochastic convergence of the quantities $\max(\htwoMC(\Yl))$ and $\max({\VltwoMC})$ is also shown in the second plot. Interestingly, one may notice that both the convergences (deterministic and stochastic) pertaining to mean and variance estimates decay in a similar manner.  
On the other hand, in the top left plot, the quantities $\max(\muMC(\UhlNO))$ and $\max(\htwoMC(\UhlNO))$ stay approximately constant at all values of $\meshlevel$. Similarly, on the right hand side, the entities   $\max(\htwoMC(\UhlNO))$ and $\max(\VtwoMC(\UhlNO))$ are also approximately constant for varying values of $\meshlevel$---which meets the assumption for consideration of standardizing factors $\Nmhat$ and $\Nvhat$ on coarse mesh $\meshlevel=0$ made in \fsecs{MLMC}{MLMCvar}. 
Furthermore, the logarithmic computational time of running one sample $\costlevel$ on each mesh level $\meshlevel$ is shown at the bottom. The corresponding values are obtained by recording the timings for the 50 considered screening samples on a 2.3GHz Intel Core i5 processor with 8GB of RAM  and taking the average. As the results show, computing becomes more expensive as the mesh refinement increases.
	
    \begin{figure}[!ht]
        \begin{subfigure}{0.49\textwidth}
            \centering
%
%
%
\begin{tikzpicture}
\footnotesize
\begin{axis}[%
width=\figWW,
height=\figHH,
at={(0\figWW,0\figHH)},
scale only axis,
xmin=1,
xmax=4,
xtick={1,2,3,4},
xticklabels={{0},{1},{2},{3}},
xlabel style={font=\color{white!15!black}},
xlabel={Mesh level, $\meshlevel$},
ymin=-6,
ymax=4,
ylabel style={font=\color{white!15!black}},
axis background/.style={fill=white},
legend style={at={(0,1.74)}, anchor=north west, legend cell align=left, align=left, draw=none}
]
\addplot [color=mycolor1, line width=1.2pt,  mark size=2.2pt, mark=square, mark options={solid, mycolor1}]
table[row sep=crcr]{%
2	-0.51430938206555\\
3	-1.77598050398772\\
4	-3.30284412562878\\
};
\addlegendentry{$\log \textrm{max}\mid\muMC(\Yl)\mid/\sqrt{\Nmhatinf}$}

\addplot [color=mycolor2, line width=1.2pt,  mark size=2.2pt, mark=o, mark options={solid, mycolor2}]
table[row sep=crcr]{%
1	2.9789582457785\\
2	3.02242958044148\\
3	3.01897245454025\\
4	3.02372425860273\\
};
\addlegendentry{$\log \textrm{max}(\muMC(\UhlNO)/\sqrt{\Nmhatinf})$}

\addplot [color=mycolor1, dashed, line width=1.2pt,  mark size=2.2pt, mark=square, mark options={solid, mycolor1}]
table[row sep=crcr]{%
	2	-3.28517026366654\\
	3	-4.29712164240019\\
	4	-5.63472487568784\\
};
\addlegendentry{$\log \textrm{max}\mid \Zl\mid/\sqrt{\Nvhatinf}$}

\addplot [color=mycolor2, dashed, line width=1.2pt,  mark size=2.2pt, mark=o, mark options={solid, mycolor2}]
table[row sep=crcr]{%
	1	-0.560833508871157\\
	2	-0.349939574271162\\
	3	-0.773627849901056\\
	4	-0.594161100272098\\
};
\addlegendentry{$\log \textrm{max}(\htwoMC(\UhlNO)/\sqrt{\Nvhatinf})$}

\end{axis}
\end{tikzpicture}%
        \end{subfigure}
        \begin{subfigure}{0.49\textwidth}
            \centering
%
%
%
\begin{tikzpicture}
\footnotesize
\begin{axis}[%
width=\figWW,
height=\figHH,
at={(0\figWW,0\figHH)},
scale only axis,
xmin=1,
xmax=4,
xtick={1,2,3,4},
xticklabels={{0},{1},{2},{3}},
xlabel style={font=\color{white!15!black}},
xlabel={Mesh level, $\meshlevel$},
ymin=-6,
ymax=2,
ylabel style={font=\color{white!15!black}},
axis background/.style={fill=white},
legend style={at={(0,1.65)}, anchor=north west, legend cell align=left, align=left, draw=none}
]
\addplot [color=mycolor1, line width=1.2pt,  mark size=2.2pt, mark=square, mark options={solid, mycolor1}]
  table[row sep=crcr]{%
2	-3.12682459181712\\
3	-3.66902109528048\\
4	-5.10520512590322\\
};
\addlegendentry{$\log \textrm{max}(\htwoMC(\YhlNl)/\Nmhatinf)$}

\addplot [color=mycolor1, dashed, line width=1.2pt,  mark size=2.2pt, mark=square, mark options={solid, mycolor1}]
table[row sep=crcr]{%
	2	-2.94176775573239\\
	3	-3.92609060247015\\
	4	-4.98525782034718\\
};
\addlegendentry{$\log {\textrm{max}(\VltwoMC/\Nvhatinf)}$}

\addplot [color=mycolor2, line width=1.2pt,  mark size=2.2pt, mark=o, mark options={solid, mycolor2}]
  table[row sep=crcr]{%
1	0\\
2	0.210890559360591\\
3	-0.212798511302271\\
4	-0.033337098803894\\
};
\addlegendentry{$\log \textrm{max}(\htwoMC(\UhlNO)/\Nmhatinf)$}

\addplot [color=mycolor2, dashed, line width=1.2pt,  mark size=2.2pt, mark=o, mark options={solid, mycolor2}]
table[row sep=crcr]{%
	1	0\\
	2	-0.211883461087118\\
	3	-0.515540210096864\\
	4	-0.362693657042124\\
};
\addlegendentry{$\log \textrm{max}(\VtwoMC(\UhlNO)/\Nvhatinf)$}

\end{axis}
\end{tikzpicture}%
        \end{subfigure} \\
        \begin{subfigure}{\textwidth}
            \centering
%
%
\begin{tikzpicture}
\footnotesize
\begin{axis}[%
width=\figWW,
height=\figHH,
at={(0\figWW,0\figHH)},
scale only axis,
xmin=1,
xmax=4,
xtick={1,2,3,4},
xticklabels={{0},{1},{2},{3}},
xlabel style={font=\color{white!15!black}},
xlabel={Mesh level, $\meshlevel$},
ymin=-2.2,
ymax=2,
ylabel style={font=\color{white!15!black}},
ylabel={$\log(\costlevel)$},
axis background/.style={fill=white},
legend style={legend cell align=left, align=left, draw=white!15!black}
]
\addplot [color=mycolor2, line width=1.2pt,  mark size=2.2pt, mark=square, mark options={solid, mycolor2}]
  table[row sep=crcr]{%
1	-1.76514636829799\\
2	-1.17665629407842\\
3	-0.323492824802114\\
4	1.64472054382001\\
};

\end{axis}
\end{tikzpicture}%
        \end{subfigure}       
        \caption{Screening results of scale-invariant MLMC mean and variance estimators}
        \label{Fig:meanscreening}
    \end{figure}


Finally, the decay rates and constants corresponding to \feqs{Eq:convNT}{Eq:convNTvar}, as well as the third equation in \feqs{eq:errormean1}{eq:errorvar1} are evaluated by determining the slopes and y-intercepts of respective logarithmic quantities in \ffig{Fig:meanscreening}, which are further summarised in \ftbl{Table:convergence_mean}. 
Clearly, all the constants are positive, and the condition $\convorder\ge\frac{1}{2}\text{min}(\convbeta,\convgamma)$ is satisfied, verifying the assumption made in \props{prop:mlmcmean}{prop:mlmcvar}. One may notice the decay rates $\convorder$ and $\convbeta$ are closer, and equal values of order $\convgamma$, for both mean and variance estimators. The constant $\cbeta$ remains close, with equal $\cgamma$ and differing $\calpha$, for both estimates. 
As $\convbeta<\convgamma$ for both mean and variance estimators, the computational complexity of MLMC estimates follows the third scenario in \feqs{Eq:costcomplex}{Eq:costcomplex2}, respectively. \rvsn{Moreover, in the considered example with a fixed number of mesh levels, the total computing cost of both MLMC estimators is dependent on the stochastic convergence order $\beta$. As the order $\beta$ for variance estimate is slightly higher than that of the mean, for a given sampling accuracy $\epssqtwo$, the estimation of MLMC variance is projected to be \rvst{more} expensive than the mean estimation.}

\begin{table}[!ht]
	\centering
	\begin{tabular}{lllllll}
		\specialrule{1pt}{1pt}{1pt}
		Statistic & $\convorder$ & $\convbeta$  & $\convgamma$ & $c_{\alpha}$  & {$c_{\beta}$}  & $c_{\gamma}$  \\ \midrule
		Mean &  2.01 &  1.43 &  \multirow{2}{*}{1.6} & 2.52 &  0.14  & \multirow{2}{*}{0.13} \\
		Variance & 1.70  & 1.47  &  & 0.13 & 0.15  &  \\
		\specialrule{1pt}{1pt}{1pt}
	\end{tabular}
	\caption{Convergence results of MLMC mean and variance estimators}
	\label{Table:convergence_mean}
\end{table}

\subsection{Performance analysis}
\ffig{performance_MLMC} shows an overview of the performance of scale-invariant \acrshort{mlmc}. 
The top plot shows the propagation of a maximum number of samples $\text{max}(\Nl)$ on each level $l$ for both MLMC mean and variance estimators, corresponding to varying normalized sampling accuracies $\epssqtwo$. One may observe that $N_l$ decreases monotonically with increasing $l$ for both estimates. 
For all of the investigated accuracies, it is obvious that the variance estimate requires more samples than the mean estimate. On the other hand, the maximum number of MC samples run on mesh \rvsb{level} $\meshlevel=3$, satisfying the given stochastic NMSEs, is tabulated in \ftbl{table:mcsample}.
\begin{figure}[h!]
	\hspace{-0.2cm}
	\begin{subfigure}{\textwidth}
		\centering
%
%
%

\begin{tikzpicture}
\footnotesize
\begin{axis}[%
width=\figWW,
height=\figHH,
at={(0\figWW,0\figHH)},
scale only axis,
xmin=1,
xmax=4,
xtick={1,2,3,4},
xticklabels={{0},{1},{2},{3}},
xlabel style={font=\color{white!15!black}},
xlabel={Mesh level, $\meshlevel$},
ymode=log,
ymin=30,
ymax=15000,
yminorticks=true,
ylabel style={font=\color{white!15!black}},
ylabel={Max. samples per level, max($\Nl$)},
axis background/.style={fill=white},
legend style={at={(1.03,1)}, anchor=north west, legend cell align=left, align=left, draw=none}
]
\addplot [color=mycolor2, line width=1.2pt, mark size=2.2pt, mark=square, mark options={solid, mycolor2}]
table[row sep=crcr]{%
1	3436\\
2	568\\
3	267\\
4	50\\
};
\addlegendentry{Mean: $\epssqtwo=6\ex{-4}$}

\addplot [color=mycolor2, line width=1.2pt, mark size=2.2pt, mark=o, mark options={solid, mycolor2}]
table[row sep=crcr]{%
1	5272\\
2	877\\
3	419\\
4	83\\
};
\addlegendentry{Mean: $\epssqtwo=4\ex{-4}$}

\addplot [color=mycolor2, line width=1.2pt, mark size=3pt, mark=x, mark options={solid, mycolor2}]
table[row sep=crcr]{%
1	10570\\
2	1782\\
3	873\\
4	162\\
};
\addlegendentry{Mean: $\epssqtwo=2\ex{-4}$}

\addplot [color=mycolor1, dashed, line width=1.2pt, mark size=2.2pt, mark=square, mark options={solid, mycolor1}]
table[row sep=crcr]{%
1	4274\\
2	967\\
3	458\\
4	100\\
};
\addlegendentry{Var: $\epssqtwo=6\ex{-4}$}

\addplot [color=mycolor1, dashed, line width=1.2pt, mark size=2.2pt, mark=o, mark options={solid, mycolor1}]
table[row sep=crcr]{%
1	6410\\
2	1514\\
3	752\\
4	137\\
};
\addlegendentry{Var: $\epssqtwo=4\ex{-4}$}

\addplot [color=mycolor1, dashed, line width=1.2pt, mark size=3pt, mark=x, mark options={solid, mycolor1}]
table[row sep=crcr]{%
1	12710\\
2	3127\\
3	1489\\
4	260\\
};
\addlegendentry{Var: $\epssqtwo=2\ex{-4}$}

\end{axis}

\begin{axis}[%
width=1.227\figWW,
height=1.227\figHH,
at={(-0.291\figWW,-0.135\figHH)},
scale only axis,
xmin=0,
xmax=1,
ymin=0,
ymax=1,
axis line style={draw=none},
ticks=none,
axis x line*=bottom,
axis y line*=left,
legend style={legend cell align=left, align=left, draw=white!15!black}
]
\end{axis}
\end{tikzpicture}%
	\end{subfigure} 
	\begin{subfigure}{0.49\textwidth}
		\centering
%
%
%
\begin{tikzpicture}
\footnotesize
\begin{axis}[%
width=\figWW,
height=\figHH,
at={(0\figWW,0\figHH)},
scale only axis,
xmin=0.0002,
xmax=0.0006,
xtick={0.0002,0.0004,0.0006},
xlabel style={font=\color{white!15!black}},
xlabel={Sampling NMSE},
ymin=0,
ymax=7,
ylabel style={font=\color{white!15!black}},
ylabel={Cost (h)},
axis background/.style={fill=white},
legend style={at={(0,1.5)}, anchor=north west, legend cell align=left, align=left, draw=none}
]
\addplot [color=mycolor2, line width=1.2pt, mark size=2.2pt, mark=o, mark options={solid, mycolor2}]
  table[row sep=crcr]{%
0.0002	1.06371973743722\\
0.0004	0.529402760700555\\
0.0006	0.337615137262778\\
};
\addlegendentry{MLMC mean in cm}

\addplot [color=mycolor2, line width=1.2pt, mark size=2.2pt, mark=square, mark options={solid, mycolor2}]
  table[row sep=crcr]{%
0.0002	6.37929343055555\\
0.0004	3.18964671527778\\
0.0006	2.12685642974722\\
};
\addlegendentry{MC mean in cm}

\addplot [color=mycolor1, dashed, line width=1.2pt, mark size=2.2pt, mark=o, mark options={solid, mycolor1}]
  table[row sep=crcr]{%
0.0002	1.54547126717389\\
0.0004	0.782690277516667\\
0.0006	0.521958852417778\\
};
\addlegendentry{MLMC variance in cm$^2$}

\addplot [color=mycolor1, dashed, line width=1.2pt, mark size=2.2pt, mark=square, mark options={solid, mycolor1}]
  table[row sep=crcr]{%
0.0002	6.37929343055555\\
0.0004	3.18964671527778\\
0.0006	2.12685642974722\\
};
\addlegendentry{MC variance in cm$^2$}

\end{axis}

\begin{axis}[%
width=1.227\figWW,
height=1.227\figHH,
at={(-0.306\figWW,-0.135\figHH)},
scale only axis,
xmin=0,
xmax=1,
ymin=0,
ymax=1,
axis line style={draw=none},
ticks=none,
axis x line*=bottom,
axis y line*=left,
legend style={legend cell align=left, align=left, draw=white!15!black}
]
\end{axis}
\end{tikzpicture}%
	\end{subfigure}       
	\begin{subfigure}{0.49\textwidth}
		\centering
%
%
%
\begin{tikzpicture}
\footnotesize
\begin{axis}[%
width=\figWW,
height=\figHH,
at={(0\figWW,0\figHH)},
scale only axis,
xmin=0.0002,
xmax=0.0006,
xtick={0.0002, 0.0004, 0.0006},
xlabel style={font=\color{white!15!black}},
xlabel={Sampling NMSE},
ymin=0.2,
ymax=1.6,
ylabel style={font=\color{white!15!black}},
ylabel={Cost (h)},
axis background/.style={fill=white},
legend style={at={(0,1.5)}, anchor=north west, legend cell align=left, align=left, draw=none}
]
\addplot [color=mycolor2,  line width=1.2pt, mark size=2.2pt, mark=o, mark options={solid, mycolor2}]
  table[row sep=crcr]{%
0.0002	1.06371973743722\\
0.0004	0.529402760700555\\
0.0006	0.337615137262778\\
};
\addlegendentry{MLMC mean in cm}

\addplot [color=mycolor1, dashdotted,  line width=1.2pt, mark size=2.2pt, mark=o, mark options={solid, mycolor1}]
  table[row sep=crcr]{%
0.0002	1.54547126717389\\
0.0004	0.782690277516667\\
0.0006	0.521958852417778\\
};
\addlegendentry{MLMC variance in cm$^2$}

\addplot [color=mycolor2, line width=1.2pt, mark size=2.2pt, mark=square, mark options={solid, mycolor2}]
  table[row sep=crcr]{%
0.0002	1.06371973743722\\
0.0004	0.529402760700555\\
0.0006	0.337615137262778\\
};
\addlegendentry{MLMC mean in m}

\addplot [color=mycolor1, dashed, line width=1.2pt, mark size=2.2pt, mark=square, mark options={solid, mycolor1}]
  table[row sep=crcr]{%
0.0002	1.54547126717389\\
0.0004	0.782690277516667\\
0.0006	0.521958852417778\\
};
\addlegendentry{MLMC variance in m$^2$}

\end{axis}
\end{tikzpicture}%
	\end{subfigure}   
	\caption{Performance of scale-invariant MLMC mean and variance estimators}
	\label{performance_MLMC}
\end{figure}
\begin{table}[!ht]
	\centering
	\begin{tabular}{ll}
		\specialrule{1pt}{1pt}{1pt}
		NMSE, $\epssqtwo$ & Number of samples, $\samplesize$   \\ \midrule
		6$\ex{-4}$ &  1667  \\
		4$\ex{-4}$ & 2500  \\
		2$\ex{-4}$ & 5000\\
		\specialrule{1pt}{1pt}{1pt}
	\end{tabular}
	\caption{Number of MC samples on level $\meshlevel = 3$ with varying stochastic accuracies}
	\label{table:mcsample}
\end{table}

The bottom left plot compares the overall cost of MLMC and MC mean and variance estimation to given normalized sampling errors.
Certainly, MLMC estimates have a faster convergence rate than the MC approach.
Furthermore, the MC cost of mean and variance exhibits non-asymptotic behaviour. However, the cost of MLMC estimates differs, with the variance estimate being relatively more expensive than the mean---as stated in \fsec{MLMCvar}.
As scale-invariant error estimates are used in this investigation, it is possible to conclude that the cost of the MLMC mean and variance estimator is asymptotic.

The MLMC estimators for the mean and variance are observed to be approximately 6 and 4 times more efficient, respectively, than their standard MC counterparts across all levels of normalized sampling accuracy.
In other words, approximate cost savings of 83\% and 75\% are reported by the mean and variance of MLMC estimates, respectively, as compared to the MC costs. 
The bottom-right graphic also shows the scale invariance element of the MLMC convergence for an individual statistic in two distinct units (cm and m). Both the MLMC mean (in cm and m) and variance (in cm$^2$ and m$^2$) convergence remain unchanged.

Finally, the summary of the given and estimated maximum value of stochastic NMSEs of MLMC mean and variance estimators is listed in \ftbl{Tbl:errors}. It can be seen that the achieved maximum NMSEs of mean and variance are within the given limits of $\epssqtwo$. Furthermore, the maximum absolute sampling accuracies, {a product of normalizing constant and maximum sampling NMSE,}
are also presented in the table. For a given NMSE value, the absolute MSE of variance is much smaller in magnitude as compared to the mean; however, both belong to different scales. This emphasizes the need for normalized error estimates to ensure easier interpretation of performance between the MLMC estimators.
\begin{table}[!ht]
	\centering
	\begin{subtable}{0.48\textwidth}
		\centering
		\begin{tabular}{lll}
			\specialrule{1pt}{1pt}{1pt}	
			{\begin{tabular}[c]{@{}l@{}}Given \\ NMSE, \\ $\epssqtwo$ \end{tabular}} & {\begin{tabular}[c]{@{}l@{}}Estimated \\ NMSE \end{tabular}} & {\begin{tabular}[c]{@{}l@{}}Estimated \\ MSE \\ ($\text{cm}^2$)\end{tabular}} \\
			\midrule
			2$\ex{-4}$                                                                             & 1.97$\ex{-4}$                                                                     & 9.459$\ex{-10}$                                                                                                                                                      \\
			4$\ex{-4}$                                                                             & 3.99$\ex{-4}$                                                                     & 1.889$\ex{-9}$                                                                                                                                                      \\
			6$\ex{-4}$                                                                             & 5.98$\ex{-4}$                                                                     & 2.859$\ex{-9}$                                                                          \\
			\specialrule{1pt}{1pt}{1pt}	
		\end{tabular}
		\caption{Mean}
	\end{subtable}
	\begin{subtable}{0.48\textwidth}
		\centering
		\begin{tabular}{lll}
			\specialrule{1pt}{1pt}{1pt}	
			{\begin{tabular}[c]{@{}l@{}}Given\\ NMSE, \\ $\epssqtwo$ \end{tabular}} & {\begin{tabular}[c]{@{}l@{}}Estimated\\ NMSE \end{tabular}} & {\begin{tabular}[c]{@{}l@{}}Estimated \\ MSE \\ ($\text{cm}^4$)\end{tabular}} \\
			\midrule
			2$\ex{-4}$                                                                           & 1.98$\ex{-4}$                                                                   & 9.387$\ex{-15}$                                                                                  \\
			4$\ex{-4}$                                                                           & 3.89$\ex{-4}$                                                                   & 1.856$\ex{-14}$                                                                                 \\
			6$\ex{-4}$                                                                           & 5.83$\ex{-4}$                                                                   & 2.789$\ex{-14}$  \\
			\specialrule{1pt}{1pt}{1pt}	
		\end{tabular}
		\caption{Variance}
	\end{subtable}
	\caption{Given and estimated sampling NMSEs and their corresponding absolute MSEs of MLMC mean and variance estimates}
	\label{Tbl:errors}
\end{table}

\begin{figure}[h!]
	\centering
	\begin{subfigure}{0.33\textwidth}
		\centering
		\includegraphics[width=0.75\textwidth]{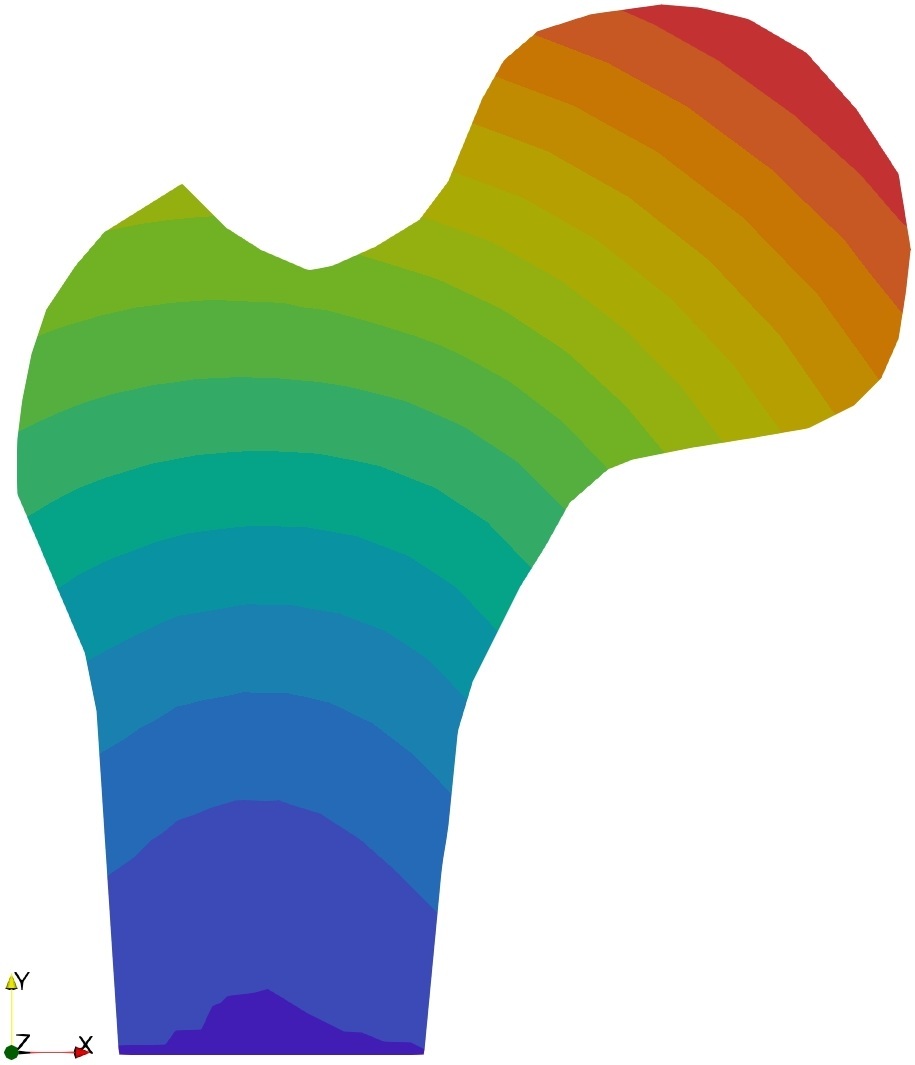}
	\end{subfigure} 
	\begin{subfigure}{0.32\textwidth}
		\centering
		\includegraphics[width=0.75\textwidth]{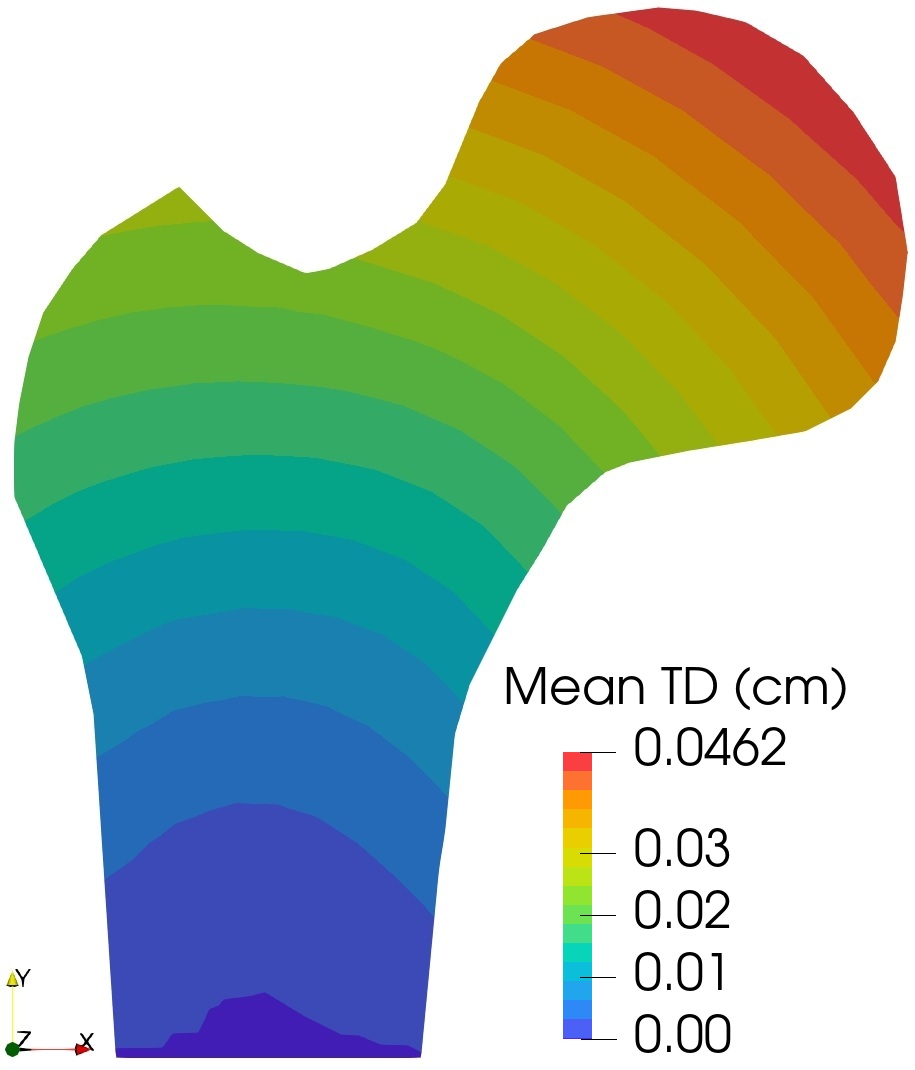}
	\end{subfigure}
	\begin{subfigure}{0.33\textwidth}
		\centering
		\includegraphics[width=0.74\textwidth]{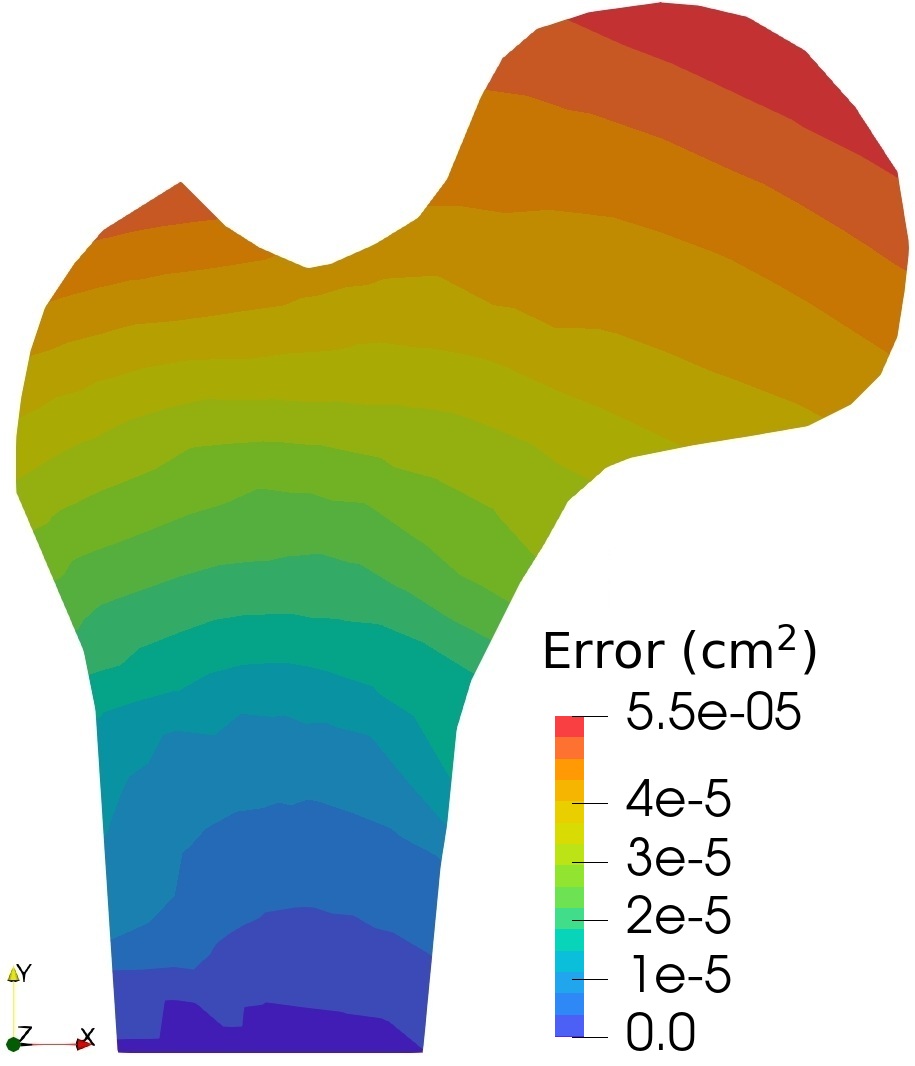}
	\end{subfigure}
	\caption{\acrshort{mlmc} (left) and \acrshort{mc} (middle) mean estimate of total displacement \rvsg{(with identical scale)}, \rvsn{along with absolute error on the right}}
	\label{Fig:mean}
	\begin{subfigure}{0.33\textwidth}
		\centering
		\includegraphics[width=0.75\textwidth]{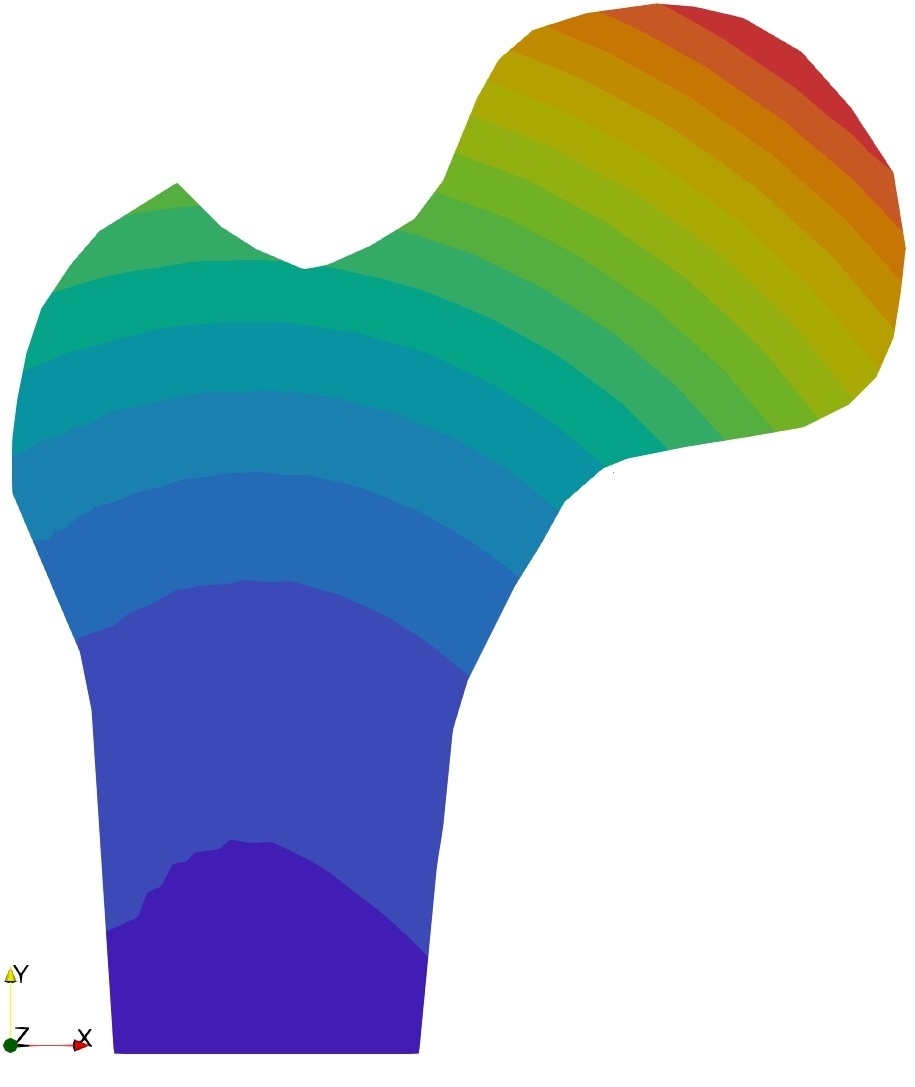}
	\end{subfigure} 
	\begin{subfigure}{0.32\textwidth}
		\centering
		\includegraphics[width=0.75\textwidth]{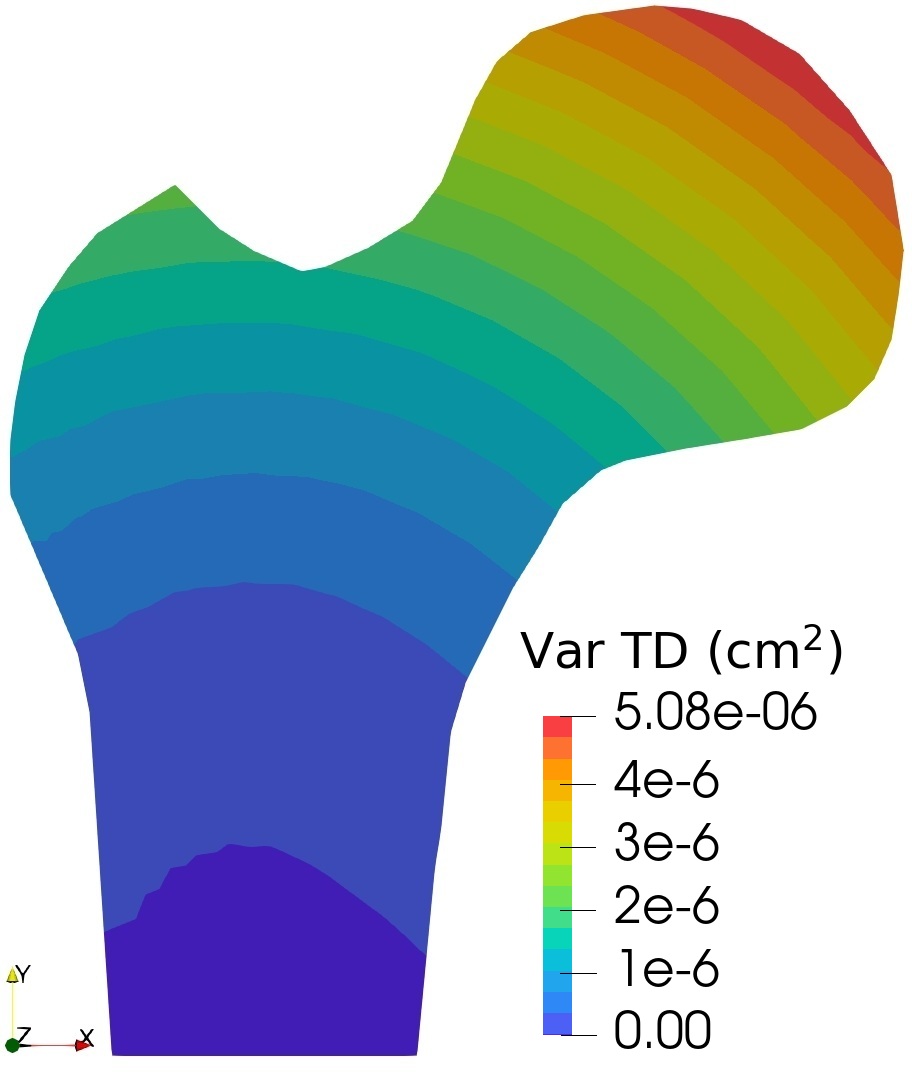}
	\end{subfigure}
	\begin{subfigure}{0.33\textwidth}
		\centering
		\includegraphics[width=0.74\textwidth]{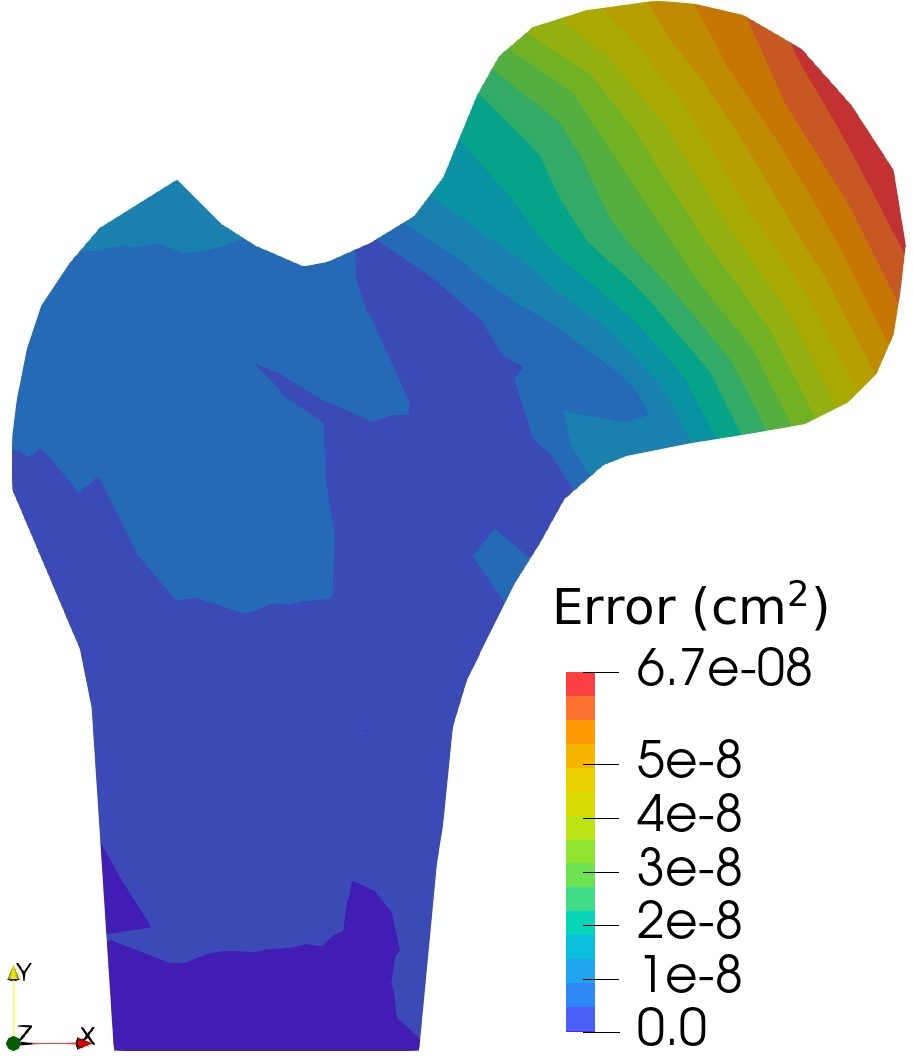}
	\end{subfigure}
	\caption{\acrshort{mlmc} (left) and \acrshort{mc} (middle) variance estimate of total displacement \rvsg{(with identical scale)}, \rvsn{along with absolute error on the right}}
	\label{Fig:variance}
\end{figure}

\ffigs{Fig:mean}{Fig:variance} compare the mean and variance estimates of total displacement (TD) of the femoral bone between MC and MLMC. \rvsg{For easier interpretation, MC and MLMC results are displayed on a single scale for both mean and variance estimates.}
In general, the results are displayed for the normalized mean-squared accuracy of $2\ex{-4}$. Note that, as the displacement values are determined on the finest mesh $\LL=3$ at the common nodes corresponding to the coarse mesh $l=0$, the contour plots are mapped directly on this coarse mesh.
In \ffig{Fig:mean}, a maximum mean value of approximately \rvst{0.0462 cm} can be seen in the region of the pressure load applied. Further, \ffig{Fig:variance} shows the influence of material uncertainty on the displacement $\uhLNLt$ due to the random material model $\CRF$. A maximum variance is also witnessed at the top region of the bone.
Furthermore, the right-most plots in both figures clearly demonstrate that the absolute error between the MLMC and standard MC methods is minimal, thereby highlighting the high level of accuracy attained by the MLMC estimators.

\section{Conclusion}\label{Conclusion}

\rvsn{We present novel scale-invariant approaches for estimating the mean and variance of a quantity of interest (QoI) using the multilevel Monte Carlo (MLMC) method. These methods are based on the derivation of normalized mean square error (NMSE) estimates for the classical Monte Carlo (MC) mean and variance estimation.
The proposed relative errors are statistically defined {using h-statistics}, with chosen normalizing factors that are finite and unbiased.
Unlike traditional normalization approaches that rely on the squared value of the estimator, which fail to achieve full scale-invariance, the proposed NMSEs are invariant under any linear transformation (scaling and translation) of QoI, and remain robust to changes in its distributional characteristics.
Such a standardized formulation reduces interpretational ambiguity and enables a dimensionless assessment of statistical accuracy and computational efficiency across different estimators and scales. }


The proposed scale-invariant MC and MLMC methods are tested on a two-dimensional simulation of a human femur modelled as an uncertain linear elastic constitutive law. As bone tissue is a highly heterogeneous and anisotropic material, and that its precise elastic symmetry class is typically unknown or uncertain, the material’s elasticity tensor is modelled as a matrix-valued random field. 
This modelling framework captures both the spatial variability of material properties and random anisotropy by prescribing an elastic symmetry in the mean (e.g., orthotropic) and allowing for triclinic symmetry in individual realizations.
The proposed methods then propagate these uncertainties to estimate statistics like the mean and variance of the stochastic total displacement field.



Through normalized error estimates, we \rvso{compare} the computational efficiencies of MLMC and MC estimators for both mean and variance. MLMC \rvso{significantly} outperforms MC in terms of computational cost for both estimates. However, the variance estimate in MLMC requires a higher number of samples, making it more computationally expensive than the mean estimate. The complete normalization of sampling error reveals non-asymptotic behaviour in the MC cost between mean and variance, while MLMC costs exhibit an asymptotic relationship. Additionally, the difference in accuracy of mean and variance estimates between the MLMC and MC methods is found to be very small.

\appendix

\section*{Appendix}

\begin{appendices}
	
\section{h-statistics: unbiased estimation of central moments} \label{appendix:hstat}

\rvsn{The h-statistics, denoted by $\hrMC$, serves as an unbiased estimator of the central population moments $\mu_\p(\uh)$ in \feq{math_Exp_h} \cite{dwyer1937, Colin_Rose:2002}. The \textquoteleft MC\textquoteright\ superscript of $\widehat{\mathrm{h}}_\p$ signifies that these statistics are obtained through random Monte Carlo-based sampling.  For example, the first two h-statistics can be defined as
	\begin{equation}\label{Eq:h-stat}
		\begin{aligned}
			\mu_1 \approx \widehat{\mathrm{h}}_1^{\text{MC}}  &= 0, \\
			\mu_2 \approx \htwoMC &= \frac{N\s_2-\s_1^2}{N(N-1)}, \\
		\end{aligned}
	\end{equation}
	in which $\s_1$ and $\s_2$ are the power sums, given as 
	\begin{equation}
		\s_a(\uh) := \sum_{k=1}^{\N}\uh(x, \event_k)^a,
	\end{equation} for $a\in\mathbb{Z}_{\geq 0}$. Using the power series representation of h-statistics, the authors in \cite{Colin_Rose:2002} developed a \textit{Mathematica}-based package called \textit{mathStatica}, which efficiently generates h-statistics for any value of $\p$.}

\rvsn{What distinguishes h-statistics from other unbiased estimators are its properties, such as \cite{halmos1946}: 
	\begin{enumerate}
		\item \textbf{Unbiasedness}: The expectation of the h-statistic equals the corresponding population central moment, i.e., $\mathbb{E}(\hrMC) = \mu_\p$.
		\item \textbf{Symmetry}: Among all unbiased estimators of $\mu_\p$, $\hrMC$ is the only one that exhibits symmetry\footnote{A function or an estimate is considered symmetric when it is not influenced by the order in which observations are considered.}.
		\item \textbf{Minimum variance}: Out of all unbiased estimators of $\mu_\p$, $\hrMC$ stands out for its minimal variance, denoted as $\var(\hrMC)$.
\end{enumerate}}

\section{Unbiased estimation of $\mu_4$ and $\mu_2^2$}\label{appendix:unbiasMCvar}
The unbiased approximation of the fourth central moment $\mu_4(\uh)$ is computed by the fourth h-statistic \cite{dwyer1937}:
\begin{multline}\label{Eq:mu_4}
	\mu_4(\uh) \approx \text{h}_4^{\text{MC}}(\uh) := \frac{1}{N(N - 1)(N - 2)(N - 3)} \bigg( -3\s_1^4 + 6N\s_1^2\s_2 \\ + (9-6N)\s_2^2 + (-4N^2+8N-12)\s_1\s_3
	+ (N^3 - 2N^2 + 3N)\s_4 \bigg),
\end{multline}
whereas the unbiased estimator of $\mu_2(\uh)^2$ is the polyache $\text{h}_{\left\lbrace 2,2\right\rbrace}^{\text{MC}}$ according to \cite{tracy1974,Colin_Rose:2002}, given as
\begin{multline}\label{Eq:mu_2_2}
	\mu_2(\uh)^2 \approx \text{h}_{\left\lbrace 2,2\right\rbrace}^{\text{MC}}(\uh) := \frac{1}{N(N - 1)(N - 2)(N - 3)}\bigg(\s_1^4 - 2N\s_1^2\s_2  \\+ (N^2 - 3N + 3)\s_2^2 + (4N - 4)\s_1\s_3+ (- N^2 + N)\s_4 \bigg).
\end{multline}
It is to be noted that the estimator of $\mu_2^2$ by the square of the second h-statistic, i.e., $(\htwoMC)^2$,  leads to a biased estimate, whereas only the polyache $\text{h}_{\left\lbrace 2,2\right\rbrace}^{\text{MC}}$, which is unbiased, remains the preferred way.
%

\section{\rvst{Variance of normalizing factor $\Nvhat:=\VtwoMC(\uh)$}}\label{appendix:varV2}

\rvst{One may define the variance of the normalizing factor $\Nvhat := \VtwoMC(\uh)$ using the statistical package \textit{mathStatica} \cite{Colin_Rose:2002} as
\begin{multline}\label{Eq:varV2}
	\var(\VtwoMC(\uh)) = \frac{72 \mu _2^4 \left(N^2-6 N+12\right)}{(N-3) (N-2) (N-1) N}+\frac{16 \mu _3^2 \mu _2 \left(N^2-4 N+13\right)}{(N-2) (N-1) N}\\-\frac{24 \mu _4 \mu _2^2 (4 N-11)}{(N-2) (N-1) N}+\frac{16 \mu _6 \mu _2}{(N-1) N}+\frac{\mu _8}{N}-\frac{8 \mu _3 \mu _5}{N}-\frac{\mu _4^2 (N-17)}{(N-1) N}.
\end{multline}
	Here, the $p$-th central moment is denoted by $\mu_p \equiv \mu_p(\uh)$. For brevity, the previous expression is rewritten as
\begin{equation}\label{Eq:simplvarV2}
	\var(\VtwoMC(\uh)) = \frac{\mathcal{V}_2}{N},
\end{equation}
where,
\begin{align}
	\mathcal{V}_2 
	&= \frac{72 \mu _2^4 \left(N^2 - 6N + 12\right)}{(N-3)(N-2)(N-1)}
	+ \frac{16 \mu _3^2 \mu _2 \left(N^2 - 4N + 13\right)}{(N-2)(N-1)} \nonumber\\
	&\quad - \frac{24 \mu _4 \mu _2^2 (4N - 11)}{(N-2)(N-1)}
	+ \frac{16 \mu _6 \mu _2}{N-1} 
	+ \mu _8 - 8 \mu _3 \mu _5
	- \frac{\mu _4^2 (N - 17)}{N - 1}.
\end{align}
\feq{Eq:simplvarV2} indicates that the statistical error of the estimator $\VtwoMC(\uh)$ converges at a rate of $\mathcal{O}(N^{-1})$.}

\section{Unbiased estimation of $\Vltwo$}\label{appendix:unbiasMLvar}
Following the notion of power sum $\s_{a}$ defined in \app{appendix:hstat}, here we introduce the bivariate power sum \cite{Colin_Rose:2002,krumscheid_quantifying_2020}
\begin{equation}
	\s_{a,b} := \sum_{i=1}^{\N} {X_{\h_l}^+(\event_i)}^a {X_{\h_l}^-(\event_i)}^b,
\end{equation}	
where 
$ X_{\h_l}^+(\event_i)_{i=1,...,\Nl} := X_{\h_l,\Nl}^+ = \uhlNl + \uhllNl $
and
$ X_{\h_l}^-(\event_i)_{i=1,...,\Nl} := X_{\h_l,\Nl}^- = \uhlNl - \uhllNl $. 
Subsequently, the unbiased MC estimation of the quantity $\Vltwo$ reads:
\begin{equation}
	\begin{split}
		{\VltwoMC} = \frac{1}{(\Nl-3)(\Nl-2)(\Nl-1)^2 \Nl}\Big( \Nl\big( ( -\Nl^2+\Nl+2){(\s_{1,1}})^2  \\
		+(\Nl-1)^2 (\Nl \s_{2,2}-2\s_{1,0}\s_{1,2}) + (\Nl-1)\s_{0,2}(({\s_{1,0}})^2-\s_{2,0})\big) \\
		+ ({\s_{0,1}})^2\big((6-4\Nl)({\s_{1,0}})^2 + (\Nl-1)\Nl \s_{2,0}\big)- \\ 
		2\Nl \s_{0,1}\big((\Nl-1)^2 \s_{2,1}+(5-3\Nl)\s_{1,0}\s_{1,1}\big) \Big).
	\end{split}	
\end{equation}		

\end{appendices}

\section*{Acknowledgment}
The authors gratefully acknowledge the financial support of the German Research Foundation (DFG) within the DFG Priority Program SPP 1748 “Reliable Simulation Techniques in Solid Mechanics”.


\bibliographystyle{abbrv-ssk-doi.bst}
\bibliography{./misc/references.bib}

\end{document}